\documentclass[a4paper,11pt,reqno]{amsart}

\usepackage{calc}

\ifx\pdftexversion\undefined
  \usepackage[dvips]{graphicx}
  \DeclareGraphicsExtensions{.eps,.ps,.mps}
\else
  \usepackage[pdftex]{graphicx}
  \DeclareGraphicsExtensions{.pdf,.png,.mps}
\fi

\usepackage{hyperref}
\hypersetup{
    pdftoolbar=true,     
    pdfstartview=XYZ,    
    bookmarksopen=true,  
    linktocpage=false,   
    colorlinks=false,    
    pdfcreator={pdflatex},
    pdfauthor={Daniel John space \textparenleft Bochum \textparenright},
    pdftitle={Symmetrization Procedures for the Isoperimetric Problem 
      in Symmetric Spaces of Noncompact Type},
    pdfsubject={Differential Geometry \space
	\textparenleft Geometric Variational Problems \textparenright},
    pdfkeywords={isoperimetric problem, symmetrization procedures, 
       symmetric spaces of noncompact type},
  }

   \makeatletter        
   \renewcommand{\section}{\@startsection
        {section}
        {1}
        {0mm}
        {\baselineskip}
        {0.5\baselineskip}
        {\normalfont\large\bfseries\centering}}
   \makeatother

\theoremstyle{plain}
\newtheorem{theorem}{Theorem}[section]
\newtheorem{theorem*}{Theorem}
\newtheorem{proposition}[theorem]{Proposition}
\newtheorem{lemma}[theorem]{Lemma}
\newtheorem{corollary}[theorem]{Corollary}

\newtheorem{definition}[theorem]{Definition}
\newtheorem{definition*}{Definition}
  
\theoremstyle{remark}
\newtheorem{remark}[theorem]{Remark}

\numberwithin{equation}{section}

\newcommand{\abs}[1]{\lvert#1\rvert}			
\newcommand{\norm}[1]{\lVert#1\rVert}			
\newcommand{\scp}[2]{\langle#1,#2\rangle}		

\newcommand{\const}{\mathrm{const}}                     
\newcommand{\esssup}{\operatorname*{ess\,sup}}          
\newcommand{\supp}{\operatorname*{supp}}                
\newcommand{\osc}{\operatorname*{osc}}                  

\DeclareMathOperator{\nab}{\nabla}			

\newcommand{\ddt}{\tfrac{d}{dt}}
\DeclareMathOperator{\Div}{div}                         
\DeclareMathOperator{\grad}{grad}			
\DeclareMathOperator{\Hess}{Hess}			

\newcommand{\area}{\operatorname*{area}}                
\newcommand{\vol}{\operatorname*{vol}}                  
\newcommand{\dist}{\operatorname*{dist}}                
\newcommand{\Ric}{\operatorname*{Ric}}                  
\newcommand{\dvol}{\mathrm{dvol}}                       
\newcommand{\Hm}{\operatorname{\mathcal{H}}}            
\newcommand{\dHm}{\mathrm{d}{\mathcal{H}}}              
\newcommand{\Per}{\operatorname*{Per}}                  
\newcommand{\Hor}{\operatorname*{Hor}}                  

\newcommand{\bbR}{\mathbb{R}}                           
\newcommand{\bbC}{\mathbb{C}}                           
\newcommand{\Sph}{\mathbb{S}}                           
\newcommand{\Hyp}{\mathbb{H}}                           
\newcommand{\BV}{\mathrm{BV}}                           
\newcommand{\loc}{\mathrm{loc}}                         
\newcommand{\LP}{\mathcal{L}}                           

\newcommand{\Isom}{\mathsf{Isom}}                       
\newcommand{\GG}{\mathsf{G}}                            
\newcommand{\GK}{\mathsf{K}}                            

\newcommand{\hM}{\widehat{M}}
\newcommand{\hN}{\widehat{N}}

\newcommand{\hc}{\widehat{c}}
\newcommand{\hg}{\widehat{g}}
\newcommand{\hOmega}{\widehat{\Omega}}
\newcommand{\hvarphi}{\widehat{\varphi}}

\newcommand{\AF}{\operatorname{\mathcal{F}}}     
\newcommand{\df}{\widetilde{a}}                  

\newlength{\awidth}
\newcommand{\ddf}{\widetilde{
                  \parbox[b][\height-0.25ex][b]{\awidth}{$\widetilde{a}$}}
                  \hspace{0.01em}} 

\newlength{\bwidth}
\newcommand{\ddb}{\widetilde{
                  \parbox[b][\height-0.25ex][b]{\bwidth}{$\widetilde{b}$}}}

\begin{document}


\title[Symmetrization Procedures for the Isoperimetric Problem]
      {\bf Symmetrization Procedures\\
       for the Isoperimetric Problem\\
       in Symmetric Spaces of Noncompact Type}
\author{Daniel John}
\address{Ruhr-Universit\"at Bochum\\
	 Fakult\"at f\"ur Mathematik\\
	 Universit\"atsstr.~150\\
	 D--44780 Bochum,
	 Germany}
\email{john@math.ruhr-uni-bochum.de}
\thanks{The author was supported by the DFG Schwerpunkt 1154 ``Globale
      Differential\-geometrie''. Special thanks are due to 
      Prof.~Dr.~Uwe Abresch for his support during the work on this thesis.}
\date{\today}


\begin{abstract}
We establish a new symmetrization procedure for the
iso\-perimetric problem in symmetric spaces of noncompact type. 
This symmetrization generalizes the well known Steiner
symmetrization in euclidean space.
In contrast to the classical construction the symmetrized domain is obtained
by solving a nonlinear elliptic equation of mean curvature type.
We conclude the paper discussing possible applications to the isoperimetric
problem in symmetric spaces of noncompact type.
\end{abstract}

\maketitle



\section*{Introduction}


In this article we consider the isoperimetric problem in symmetric spaces of
noncompact type, i.e., the problem of determining the domains minimizing
surface area among all regions with a given volume. 
As existence and partial regularity of isoperimetric solutions in these spaces
are given by geometric measure theory \cite[pp.~129]{Morgan1},
the goal here is to get some information about the shape of isoperimetric
solutions 
in these spaces. 

In the history of the isoperimetric problem symmetrization procedures have
been a very important tool. J.~Steiner (1838), H.~A.~Schwarz (1884),
and E.~Schmidt (1943) used 
symmetrization arguments to get insight into the behavior of isoperimetric
solutions in $\bbR^n$, $\Hyp^n$, and $\Sph^n$, finally proving 
the isoperimetric
property of metric balls in constant curvature spaces \cite{BuragoZalgaller1}.

Beginning in 1989 with the work of W.-T.~Hsiang and W.-Y.~Hsiang 
\cite{Hsiang1} the isoperimetric problem has been investigated in spaces like 
$\Hyp^n\times\bbR^m$, $\Hyp^n\times\Hyp^m$, $\Sph^n\times\Sph^1$,
$\bbR^n \times \Sph^1$, $\Hyp^n \times \Sph^1$, or $\Sph^n\times\bbR$
by  R.~Pedrosa, M.~Ritor\'e, and D.~John,
\cite{Pedrosa1, PedrosaRitore1, John1}. In these manifolds
the initial technical tool always is a symmetrization argument
reducing the problem to the $2$--dimensional quotient of the product space by
the isotropy group.

In some $3$--dimensional space forms, for example $\bbR P^3$, 
stability arguments have been applied successfully 
by M.~Ritor\'e and A.~Ros \cite{RitoreRos2, Ros2}.

Up to now the isoperimetric problem has been investigated only in such special
manifolds.
Techniques suitable for more general symmetric spaces are largely unknown.

The main goal of this paper is to establish a symmetrization procedure for
domains in symmetric spaces of noncompact type. So far, it is not possible to
conclude uniqueness or convexity of isoperimetric solutions by applying this
symmetrization procedure. Nevertheless, it provides some interesting insights
into the qualitative behavior of isoperimetric solutions.


\section{Main Results}


One of the fundamental features of symmetric spaces is the existence of 
special $1$--parameter groups $\tau_t$ of isometries called transvections.
Our main idea is to use these $1$--parameter groups in order to 
establish a symmetrization procedure.

\begin{figure}
\includegraphics[width=0.85\textwidth]{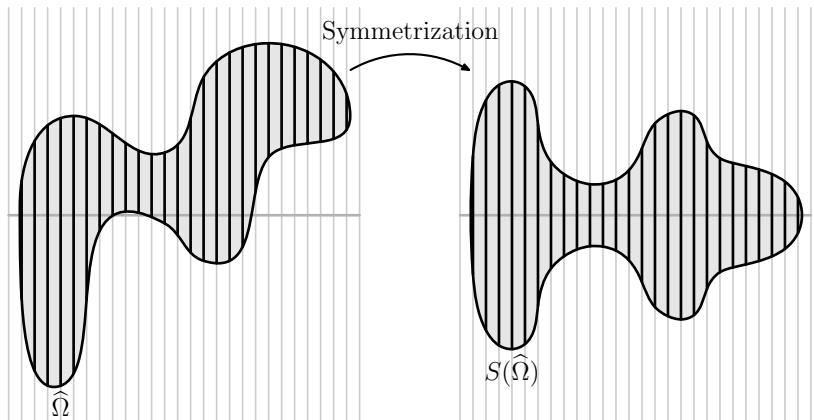}
\caption{Symmetrization using transvections}
\label{3}
\end{figure}

\begin{definition*}[\bf Symmetrization]
\label{59}
Let $\hM^n$ be a symmetric space of noncompact type, 
$\hOmega \subset \hM^n$ a given domain, and $\tau$ a transvection.
Symmetrization of $\hOmega$ with respect to $\tau$ 
is defined to be the following:
Determine a domain $S(\hOmega)$ minimizing surface area among all
volume preserving deformations of $\hOmega$ obtained by moving the line
segments
$\tau_{\bbR}(x) \cap \hOmega$, $x \in \hM^n$, along the orbits of $\tau$,
compare Figure \ref{3}.
\end{definition*}

This obviously is a generalization of the well known Steiner symmetrization,  
since the orbits of transvections in $\bbR^n$ are just parallel lines.
In euclidean space existence, uniqueness, and regularity 
properties of the symmetrized domain $S(\hOmega)$ are immediate consequences 
of Minkowski's inequality.
Establishing these properties for $S(\hOmega)$ in general symmetric spaces
of noncompact type is much more involved and one of the main issues of this 
paper.

We will mainly consider domains
$\hOmega \subset \subset \hM^n$ for which $\tau_{\bbR}(x)\cap \hOmega$
consists of a connected line segment for any $x \in \hM^n$.
These domains admit a simple representation in terms of 
the orbit space $M^{n-1}=\smash[b]\hM^n/\tau$,
a section $\sigma:M^{n-1} \to \hM^n$, 
and appropriate functions $u,h:\Omega \to \bbR$, $h \geq 0$,
where $\Omega:=\hOmega / \tau$ denotes the quotient domain:
\begin{equation*}
  \hOmega = \{\tau_t(\sigma(x)) \;|\; u(x)-h(x) \leq t 
 \leq u(x)+h(x), x \in \Omega\}.
\end{equation*}

With these notations the main theorem of the present article can 
be stated as follows.

\begin{theorem*}[\bf Symmetrization]
\label{54}
Let $\hM^n = \GG / \GK$ be a symmetric space of noncompact type.
Consider a regular domain $\hOmega \subset\subset \hM^n$ and a 
transvection $\tau$ such that the following holds:
\begin{enumerate}
\item $\tau_{\bbR}(x)\cap\hOmega$ is connected for every $x \in \hM^n$.
\item $h:\Omega \to \bbR$,  
   $h(x):=\Hm^1(\tau_{\bbR}(\sigma(x)) \cap \hOmega)$ 
   is smooth on the subset $\Omega$ of the orbit space.
   $\Hm^1$ denotes $1$--dimensional Hausdorff measure.
\end{enumerate}
Then the symmetrization procedure of Definition \ref{59}
assigns to $\hOmega$ a unique symmetrized domain $S(\hOmega) \subset \hM^n$
of equal volume but smaller (or equal) surface area.
The boundary of 
$S(\hOmega)$ is smooth in those points $x \in \partial S(\hOmega)$ where 
$\tau_{\bbR}(x) \cap \partial \hOmega$ consists of precisely
two different points.
\end{theorem*}

This theorem can easily be extended to a more general class of domains 
$\hOmega \subset \hM^n$ for
which assumption (1) does not
hold for every $x \in \hM^n$. Compare Section \ref{4} for further details.

Regularity of $\partial S(\hOmega)$ is only investigated for those
points described in the theorem but may also hold for other points. 

\begin{remark}
The word symmetrization is used since the construction resembles 
Steiner symmetrization.
However, in contrast to the case of constant curvature spaces,
$S(\hOmega)$ is not necessarily invariant under a larger class of isometries
than $\smash{\hOmega}$.
Moreover, in the general case the hypersurface defined by the midpoints
of the segments $\tau_{\bbR}(x) \cap S(\hOmega)$ does depend on the height
function $h$. It is not a priori given,
whereas for domains in constant curvature spaces it is just a hyperplane.
\end{remark}

Once we have established this symmetrization procedure we can apply it to
the isoperimetric problem in $\hM^n$.
In analogy to Steiner symmetrization in $\bbR^n$ we could try to show
convexity of isoperimetric solutions in $\hM^n$. Unfortunately this 
does not work (up to now). 
The main difficulty is illustrated in Section \ref{52} where,
by direct construction, we show the following. 

\begin{theorem*}
\label{61}
Consider a transvection $\tau$ in a symmetric space of noncompact type 
$\smash{\hM^n}$ which is not a space of constant sectional curvature.
Then there exists a domain which is not convex but invariant under
symmetrization with respect to $\tau$. This domain looks like a helix
winding up in the direction of $\tau$.
\end{theorem*}

This example provides a good starting point for further investigations.
If one assumes that there exists a nonconvex isoperimetric solution in some
symmetric space of noncompact type, then our example gives some
indications how such a solution could be constructed.

However, the helix is obviously far from being an isoperimetric
solution, just think about symmetrization with respect to other transvections.
Therefore, convexity of isoperimetric solutions in symmetric spaces of
noncompact type remains a natural conjecture. Our example
shows the main difficulties we will have to deal with proving convexity of
isoperimetric domains.
\smallskip

The present paper is organized as follows:
in Section \ref{40} we shortly review the fundamental 
properties of symmetric spaces and transvections and use them to describe 
the symmetrization procedure in more detail.
It turns out that, given $\hOmega$ as above, 
the area of $\partial \hOmega$ is computed by the functional
\begin{equation*}
  \AF(u)=\int_{\Omega} \sqrt{1+\norm{w+dh+du}^2}
                      +\sqrt{1+\norm{w-dh+du}^2}\,\dvol,
\end{equation*}
where $w$ is a $1$--form on the orbit space $M^{n-1}$.

Symmetrizing $\hOmega$ with respect to $\tau$
can be reduced to finding a function $u_0$ minimizing $\AF$.
Using $u_0$ the symmetrized set $S(\hOmega)$ is given by 
\begin{equation*}
S(\hOmega) = \{\tau_t(\sigma(x)) \;|\; u_0(x)-h(x) \leq t 
 \leq u_0(x)+h(x), x \in \Omega\}.
\end{equation*}

The main work now consists in establishing existence and regularity of
a minimizer $u_0$. For this purpose,
in Section \ref{60}, we investigate the analytic properties of the area 
functional in more detail.
It turns out that for the corresponding Euler-Lagrange equation  
we only have an estimate \eqref{6} providing nonuniform ellipticity.
Our existence and regularity results are based on this estimate
and involve ideas due to Ladyzhenskaya and Ural'tseva \cite{LU1} as well as 
Giaquinta, Modica, and Soucek \cite{GMS}.

The main issue of Section \ref{37} is to establish an a priori gradient
estimate for minimizers of $\AF$.
Finally, the work of Sections \ref{40}, \ref{60},
and \ref{37} is summarized in Theorem \ref{58} which
immediately implies Theorem \ref{54}.
\smallskip

In Section \ref{44} we prove Theorem \ref{61}.
Furthermore, the $1$--form $w$ involved in the definition of
$\AF$ is calculated explicitly for the
case of complex hyperbolic space. Resorting to the properties of $w$
in more detail should be important for future investigations of the
isoperimetric problem in these spaces.


\section{Symmetrization in symmetric spaces}
\label{40}



\subsection{Symmetric spaces and transvections}
\label{36}


Let $(\hM^n,\hg)$ be a Riemannian manifold of dimension $n$.
$(\hM^n,\hg)$ is called locally symmetric, if for every
point $x \in \hM^n$ there exists a neighbourhood $U_x$ of $x$ and a
geodesic isometry $I_x:U_x \to U_x$ such that
\begin{equation*}
  I_x(x)  =  x \quad \text{and} \quad
  {d\,I_x}_{|x}  =  -\mathrm{id}: T_x\,U_x \to T_x\,U_x. 
\end{equation*}
$(\hM^n,\hg)$ is called (globally) symmetric 
if $U_x=\hM^n$ for every $x \in \hM^n$.

It is well known that 
for each locally symmetric space $\hN^n$ there exists a simply connected
symmetric space $\hM^n$ and a group $\Gamma$ operating on 
$\hM^n$ discretly, without fixed points, and isometrically,
such that $\hN^n=\hM^n/\Gamma$.

Isoperimetric domains in locally symmetric spaces
$\hM^n/\Gamma$ depend strongly on the group $\Gamma$.
Therefore the isoperimetric problem in these spaces is mainly unsolved 
even in the easiest case of $\hM=\bbR^3$
and $\Gamma$ one of the cristallographic groups \cite{Ros1}.
In our discussion of the isoperimetric
problem we restrict our attention to the following class of symmetric spaces.

\begin{definition}
Let $\hM^n$ be a simply connected symmetric space of nonpositive
curvature. Then $\hM^n$ is said to be a symmetric space of noncompact
type, if it is not the Riemannian product of a euclidean space and another
manifold. 
\end{definition}

From now on in this work, $\hM^n$ always denotes a symmetric space 
of noncompact type. The irreducible symmetric spaces of noncompact type have
been classified by Cartan \cite[Chapter X]{Helgason1}.
\smallskip

The symmetrization procedure we will develop in this paper is based on 
special $1$--parameter groups of isometries called transvections.

\begin{definition}
Consider a normal geodesic $\gamma: \bbR \to \hM^n$. Then the 
$1$--parameter group of isometries 
$\tau_\gamma: \bbR \times \hM^n \to \hM^n$ defined
by 
\begin{equation*}
\tau_\gamma(t,x):=I_{\gamma(\frac{t}{2})} \circ I_{\gamma(0)}(x)
\end{equation*}
is called a transvection along $\gamma$. We also write $\tau(t,x)=\tau_t(x)$
for short.
The Killing field corresponding to $\tau$ is denoted by $K$.
\end{definition}

We collect some basic properties of transvections.
First of all, for all $s$, $t$ we have $\tau_t(\gamma(s))=\gamma(s+t)$ and
$\tau_{t+s} = \tau_t \circ \tau_s$.
Furthermore, $\tau$ is a differentiable map and 
$d\tau_t$ is parallel translation along $\gamma$.
For every $x \in \hM^n$ the map $t \mapsto \tau(x,t)$ is
injective because $\hM^n$ is a Hadamard
manifold. This property will be important for our surface area calculations.
It does not hold in symmetric spaces of compact type and is one of the main
reasons for restricting our attention to symmetric spaces of noncompact type.

These properties imply that on the orbit space
\begin{equation*}
  M^{n-1} := \hM^n/\tau = \left\{ \tau_\bbR(x) \;|\; x \in
  \hM^n \right\}
\end{equation*}
there exists a unique structure of a Riemannian manifold 
such that the projection $\pi: \hM^n \to M^{n-1}$,
$\pi(x) = \tau_\bbR(x)$ is a Riemannian submersion. 
\smallskip

Constructing our symmetrization procedure we will need the orbit space
$M^{n-1}$ as well as a section $\sigma: M^{n-1} \to \hM^n$. Such a section can
be obtained, for example, as follows:

Corresponding to our geodesic $\gamma: \bbR \to \hM^n$ consider the
Busemann functions $\beta_\gamma^+:\hM^n \to \bbR$
and $\beta_\gamma^-:\hM^n \to \bbR$ defined by
\begin{align*}
  \beta_\gamma^{\pm}(x)&:= \lim_{t \to \infty}(t-\dist(x,\gamma(\pm t))).
\end{align*}
Busemann functions in a symmetric space of noncompact type are known to
be $C^\infty$--differentiable. The level sets of the Busemann functions are
called horospheres. 
We now consider the function $\eta: \hM^n \to \bbR$ defined by
\begin{equation*}
  \eta = \eta_\gamma := \tfrac{1}{2}(\beta_\gamma^+ -\beta_\gamma^-).
\end{equation*}
As the isometry $\tau_t$ transfers horospheres 
$(\beta_\gamma^\pm)^{-1}(\{b\})$ into corresponding horospheres
$(\beta_\gamma^\pm)^{-1}(\{b \pm t\})$,
the function $\eta$ obviously has the following properties:
\begin{enumerate}
\item $\eta$ is $C^\infty$--differentiable,
\item $\eta(\tau(t,x))=\eta(x)+t$, and
\item $\scp{\grad \eta}{K}_{|x} = d\eta_{|x}\cdot 
  \frac{\partial}{\partial t}\tau(t,x)_{|t=0} = 1$.
\end{enumerate}
Consequently $0$ (as well as any other element of $\bbR$) is a regular value of
$\eta$ and the orbits of $\tau$ intersect the level sets of $\eta$
transversally. Summarizing we can construct the desired section.

\begin{lemma}
Let $\pi : \hM^n \to M^{n-1}=\hM^n/\tau$ be the usual projection onto the
orbit space. Then
\begin{equation*}
  \sigma:M^{n-1} \to \hM^n, \quad 
  \sigma(\pi(x)):=\eta^{-1}(\{0\}) \cap \tau(\bbR,x)
\end{equation*}
is a section and
\begin{equation*}
  \Phi: \bbR \times M^{n-1} \to \hM^n, \quad
  \Phi(t,x) := \tau(t,\sigma(x))
\end{equation*}
is a diffeomorpism.
\end{lemma}


\subsection{Surface area calculations}
\label{1}


We consider a geodesic $\gamma:\bbR \to \hM^n$ and the
corresponding $1$--parameter group $\tau$ of transvections.
Using the section $\sigma: M^{n-1} \to \hM^n$ constructed above, 
we can describe any hypersurface
in $\hM^n$ intersecting the orbits of $\tau$ transversally  
by maps
\begin{equation*}
  \psi_u: \Omega \to \hM^n, \quad \psi_u(x) := \tau(u(x),\sigma(x)),
\end{equation*}
where $\Omega \subset M^{n-1}$ and $u: \Omega \to \bbR$ is an appropriate
function. 
Using $u$ we can derive an easy formula computing
the surface area of $\psi_u(\Omega) \subset (\hM^n,\hg)$.
For this issue we need the following notation.

For a vector field $X$ on the orbit space $M^{n-1}$ we denote 
by $\Hor(X)$ the unique horizontal vector field on $\hM^n$ (with respect to
the submersion $\pi$) such that $d\pi(\Hor X) = X_{|\pi}$. 
$K$ is the Killing field corresponding to $\tau$.

\begin{definition}
The $1$--forms $w_u$ and $w$ on the subset 
$\Omega \subset M^{n-1}$ of the orbit space are defined
by requiring 
\begin{align*}
  d\psi_u \cdot X &= \Hor(X)_{|\psi_u} + w_u(X) \cdot K_{|\psi_u},\\
  d\sigma \cdot X &= \Hor(X)_{|\sigma} +w(X) \cdot K_{|\sigma},
\end{align*}
which is the splitting of 
$d\psi_u \cdot X$ and $d\sigma \cdot X$ into horizontal and vertical part.
These $1$--forms are obviously related by $w_u = w + du$.

\end{definition}

We denote the volume form of the orbit space $(M^{n-1},g)$ by
$\dvol_g$. The volume form of the hypersurface 
$\psi_u(\Omega) \subset \hM^n$
with respect to the metric induced by $\hg$ is denoted by
$\dvol_{\hg}$. Using this notation, the volume form of 
$(\Omega,\psi_u^*\hg)$ is $\psi_u^*\dvol_{\hg}$.

\begin{lemma}
\label{50}
Define $k:\Omega \to \bbR$, $k(x):=\norm{K_{|\psi_u(x)}}_{\hg}
=\norm{K_{|\sigma(x)}}_{\hg}$. Then
\begin{align*}
  \psi_u^* \dvol_{\hg} 
  &= \sqrt{1+k^2\cdot \norm{w+du}_g^2} \;\dvol_g\\
  &= \sqrt{1+k^2\cdot \norm{W+\grad u}_g^2}  \;\dvol_g,
\end{align*}
where $W$ is the vector field related to the $1$--form $w$ by
$w(X) = g(W,X)$ for all vector fields $X$ on $M^{n-1}$.
For $\Omega \subset M^{n-1}$ and $u \in C^\infty(\Omega)$ the surface area
of the set 
$\psi_u(\Omega) = \{\tau(u(x),\sigma(x)) \;|\; x \in \Omega\}
\subset \hM^n$ is computed by
\begin{align*}
  \area(u) := \area(\psi_u(\Omega))
  &= \int_\Omega \sqrt{1+k^2\cdot \norm{w+du}_g^2}  \;\dvol_g\\
  &= \int_\Omega \sqrt{1+k^2\cdot \norm{W+\grad u}_g^2}  \;\dvol_g.
\end{align*}
\end{lemma}

\begin{proof}
Writing $(\psi_u^*\hg)(X,Y)=g(X,G\cdot Y)$ for an appropriate field of
endomorphisms $G$ we have
$(\psi_u^*\dvol_{\hg})
  = \sqrt{\det{G}} \cdot \dvol_g$.
A short calculation shows that
\begin{align*}
  g(X,G\cdot Y)
  &= g(X,Y)+w_u(X) \cdot w_u(Y) \cdot 
           \hg(K_{|\psi_u},K_{|\psi_u}) \text{ and}\\
  \det G &= 1 + \norm{K_{|\psi_u}}_{\hg}^2 \cdot \norm{w_u}_g^2
    = 1 + \norm{K_{|\psi_u}}_{\hg}^2 \cdot \norm{w + du}_g^2.   
\end{align*}
\end{proof}

We establish some basic properties of the function $k$ and the $1$--form
$w$.

\begin{lemma}
\label{46}
For $k:M^{n-1} \to \bbR$, $k(x):=\norm{K_{|\sigma(x)}}_{\hg}$
we have
\begin{enumerate}
\item $k \in C^\infty(M^{n-1})$.
\item $k$ is a convex function on $M^{n-1}$ with $k \geq 1$ and 
  $k(\pi(\gamma(t))) = 1$, where $\gamma$ is the geodesic incorporated in the
  definition of $\tau = \tau_\gamma$.
\item $k$ is invariant under the isometries on the orbit space which are
  induced by those isometries of $\hM^n$ that transfer 
  $K_{|\gamma(0)}$ into $\pm K_{|\gamma(0)}$.
\end{enumerate}
\end{lemma}

\begin{proof}
(1) and (3) are clear. For (2) consider a geodesic
$c:\bbR \to M^{n-1}$ in the orbit space and a horizontal geodesic
$\hc$ in $\hM^n$ with $\pi \circ \hc = c$. Then
$k(c(t))=\norm{K_{|\sigma(c(t))}}=\norm{K_{|\hc(t)}}$.
Since $K$ is a Killing field, $t \mapsto K_{|\hc(t)}$ is a Jacobi field. 
Now ${\hM}$ has nonpositive curvature and therefore
$t \mapsto \norm{K_{|\hc(t)}}$ is a smooth, convex function.
Hence $k(x)\geq k(\pi(\gamma(t)))=1$ for all $x \in M^{n-1}$.
\end{proof}

\begin{lemma}
\label{47}
$W$ is a smooth vector field on $M^{n-1}$.
If $\hvarphi \in \Isom(\hM^n)$ is an isometry with
$d\hvarphi(K_{|\gamma(0)})=\pm K_{|\gamma(0)}$ and 
$\varphi \in \Isom(M^{n-1})$ the isometry induced by
$\hvarphi$ on the orbit space $M^{n-1}$,
then $d\varphi(W) = \pm W_{|\varphi}$.
\end{lemma}

\begin{proof}
$\hvarphi$ maps $\gamma$ to $\gamma$. Therefore the Killing field $K$
is mapped to $\pm K$. Furthermore, $\hvarphi$ leaves invariant
the foliation of $\hM^n$ by the level sets of the function
$\beta_{\gamma^+} - \beta_{\gamma^-}$. Hence the unit normal field
$\nu$ on $\sigma(M^{n-1})$ also remains invariant under $\hvarphi$ up to sign.
\end{proof}


\subsection{Construction of the symmetrization procedure}
\label{4}


To introduce our generalized symmetrization procedure we concentrate on the
following situation. Let $\hOmega \subset\subset \hM^n$ 
be a subset of the symmetric
space of noncompact type and $\tau=\tau_{\gamma}$ a transvection such that
\begin{enumerate}
\item $\tau_{\bbR}(x) \cap \hOmega$ is connected for every $x \in \hM^n$.
\item If $\Omega:=\pi(\hOmega)$ is the projection of $\hOmega$ to the orbit
      space then the function $h:\Omega \to \bbR$,
      $h(x):=\Hm^1(\tau(\bbR,\sigma(x)) \cap \hOmega)$
      is smooth. Here $\Hm^1$ denotes the $1$--dimensional Hausdorff measure.
\item $\partial \hOmega$ is smooth.
\end{enumerate}
These assumptions guarantee that $\hOmega$ can be written as
\begin{equation}
\label{12}
  \hOmega = \{\tau(t,\sigma(x)) \;|\; u(x)-h(x) \leq t 
 \leq u(x)+h(x), x \in \Omega\},
\end{equation}
where $u: \Omega \to \bbR$ is an appropriate smooth function.

According to Definition \ref{59} symmetrization of such an 
$\hOmega$ with respect to
$\tau$ amounts to finding a symmetrized set $S(\hOmega)$ having least surface
area among all deformations of $\hOmega$ of the form
\begin{equation*}
  \hOmega_v := \{\tau(t,\sigma(x)) \;|\; v(x)-h(x) \leq t 
 \leq v(x)+h(x), x \in \Omega\},
\end{equation*}
with $v: \Omega \to \bbR$ any function (in a reasonable function space).
Obviously all these domains $\hOmega_v$ have the same volume since 
$\tau$ is a $1$--parameter group of isometries. 

Therefore our surface area calculations imply that 
$S(\hOmega)=\hOmega_{u_0}$ where $u_0:\Omega \to \bbR$ minimizes
the surface area functional
\begin{equation}
\label{38}
  \AF(v) := \area(\partial \hOmega_v) = \int_\Omega f(\grad v) \,\dvol_g
\end{equation}
with $f: T\Omega \to \bbR$ defined by
\begin{align*}
  f(X) & := 
    \begin{aligned}[t]
        &\sqrt{1+k^2\norm{W_{|x}+\grad h_{|x}+X}^2_g}\\
        &+\sqrt{1+k^2\norm{W_{|x}-\grad h_{|x}+X}^2_g}
          \quad \text{for } X \in T_x\Omega.
    \end{aligned}
\end{align*}

In other words, in order to symmetrize $\hOmega$ we have to solve the
variational problem $\AF(u)=\min$. We will start investigating
the corresponding existence, uniqueness, and regularity questions 
in the next section.
But before, we want to make some remarks concerning the assumptions 
(1) -- (3).
 
\begin{remark}
Assumption (3), that is smoothness of $\partial \hOmega$, can be made 
without loss of generality. This can be justified as follows:
In geometric measure theory the isoperimetric problem is considered in the
class of sets of finite perimeter \cite{Giusti1}. These sets can be seen as a 
special case of the more general notion of currents. 
The perimeter of a measurable set 
$\hOmega \subset \hM^n$ is defined by 
\begin{equation*}
  \Per(\hOmega)
  := \inf (\lim\inf \Hm^{n-1}(\partial \hOmega_i),
\end{equation*}
where $\Hm^{n-1}$ denotes the $(n-1)$--dimensional Hausdorff measure
and the infimum is taken over all sequences of embedded $n$--dimensional
manifolds $\hOmega_i$ with smooth boundary 
$\partial \hOmega_i$ such that the characteristic functions 
$\chi_{\hOmega_i} \to \chi_{\hOmega}$ converge in $L^1$.
${\hOmega \subset \hM^n}$ is called a set of
finite perimeter, if $\Per(\hOmega)< \infty$.

Thinking about isoperimetric solutions as approximated by smooth
domains, it is (almost) sufficient to develop a symmetrization procedure for
sets $\hOmega \subset \hM^n$ with smooth boundary.

Another justification for assumption (3) 
can be given by the regularity part of geometric measure
theory: Isoperimetric solutions are smooth up to a singular set of 
codimension $\geq 7$. 
\end{remark}

\begin{remark}
\label{48}
The symmetrization procedure can be easily extended to a much larger class of
subsets of $\hM^n$ than those described above:
Just think about $\Omega$ not as a subset of the orbit space but as an
open set such that $\overline{\Omega}$ is a compact
$(n-1)$--dimensional Riemannian manifold with smooth boundary, isometrically
immersed into the orbit space $M^{n-1}$. Then consider again domains
$\hOmega \subset \hM^n$ whose boundary is described by functions $u-h$, $u+h$
on $\Omega$ the same way as above. 
Taking this point of view we can also apply our
symmetrization procedure to domains such as a thickened helix winding up in
the direction of the transvection $\tau_{\gamma}$.
However, one should be aware that it is not possible to describe every
smooth domain $\hOmega \subset \hM^n$ like this, i.e.,
using an $\Omega$ isometrically immersed 
into the orbit space $M^{n-1}$. An immediate counterexample is 
a torus $T$ in $\hM^n$ where $T \cap \gamma(\bbR)$ has two components.
\end{remark}

\begin{remark}
Another strategy to deal with the fact that an arbitrary domain 
$\hOmega \subset \hM^n$ can intersect the orbits 
$\tau_{\gamma}(\bbR,\sigma(x))$ in more than one component is to introduce
for every component functions $u_i,h_i:\Omega_i \to \bbR$ such that 
$u_i-h_i$, $u_i+h_i$ describe the corresponding part of the boundary
$\partial \hOmega$. 
Considering this account it seems appropriate to minimize the functional
$\AF$ on subsets $U \subset \Omega_i$ with respect to Dirichlet boundary
conditions on $\partial U$. For the (only minor) differences compared
to the case of free boundary values treated in this paper
we refer the reader to \cite{GMS}. 
\end{remark}

\begin{remark}
It should also be possible to generalize the symmetrization construction
to domains $\hOmega = \hOmega_u$ where $u$ is not smooth but in a more general
class of functions admitting jumps, such as $\BV(\Omega)$, the class of 
functions of bounded variation. This would be an interesting issue.
In this paper our focus is on the smooth case because we are interested in
uniqueness and regularity properties of the symmetrized set $S(\hOmega)$.
\end{remark}


\section{Analytic properties of the variational problem}
\label{60}



\subsection{The area functional}
\label{53}


Resuming the above observations we now investigate 
existence, regularity, and uniqueness of solutions 
for the minimizing problem $\AF(u)=\min$.

First of all, we have to think about the appropriate function spaces 
concerning this problem. For this purpose we shortly review the basic 
properties and notions of Sobolev spaces on Riemannian manifolds:

Given the open subset $\Omega \subset\subset M^n$
we denote by $L^p(\Omega)$ the space of measurable functions 
$f$ on $\Omega$ for which 
$\int_\Omega \abs{f}^p \,\dvol_g < \infty$, where $1 \leq p < \infty$. 
For a vector field $X$ on $\Omega$ we define the norms
$\norm{X}_{\LP^p}:=\left( \int_\Omega \norm{X}^p 
                            \,\dvol_g \right)^{1/p}$,
$1 \leq p < \infty$.
As usual, $\LP^p(\Omega)$ is the space of measurable vector fields $X$ 
with $\norm{X}_{\LP^p} < \infty$.

\begin{definition}
Given a function $f \in L^1_{\loc}(\Omega)$ we say that 
$Y \in \LP^1_{\loc}(\Omega)$ is a weak derivative of $f$ if
\begin{equation*}
\int_{\Omega} \scp{X}{Y} \,\dvol_g = - \int_{\Omega} f \cdot \Div X \,\dvol_g
\end{equation*}
for all $C^1$--vector fields $X$ with compact support in $\Omega$. 
If such an $Y$ exists, it is unique and we write $\grad f := Y$.
\end{definition}

\begin{definition}
The Sobolev space $W^{1,p}(\Omega)$ consists of all those 
functions $f \in L^p(\Omega)$ for which the weak derivative exists and 
$\grad f \in \LP^p(\Omega)$. 
For $f\in W^{1,p}(\Omega)$ we define its norm to be
$\norm{f}_{W^{1,p}} := (\norm{f}_{L^p}^p 
                       + \norm{\grad f}_{\LP^p}^p)^{1/p}$.
\end{definition}

As in the euclidean case, 
the usual Sobolev inequalities and embedding theorems also hold on
manifolds \cite{Aubin1}. 
Isoperimetric inequalities are closely related
to Sobolev inequalities, more precisely to the optimal constant in the Sobolev
inequality, compare \cite[p.~39]{Aubin1}.
This has been of fundamental importance in the history of the Yamabe problem
\cite[p.~153]{Aubin1}.
\smallskip

Concerning the appropriate function spaces for our area functional $\AF$ we 
consider the following lemma.

\begin{lemma}
\label{10}
Let $\scp{\cdot}{\cdot}$ be an arbitrary scalar product 
on $\bbR^m$, $\norm{\cdot}$ the
corresponding norm. Then for every $a, b \in \bbR^m$:
\begin{align*}
  \norm{a} &\leq \sqrt{1+\norm{a}^2}
  \leq \sqrt{2} \cdot \sqrt{1+\norm{a}^2+\norm{b}^2}\\
  &\leq \sqrt{1+\norm{a+b}^2} + \sqrt{1+\norm{a-b}^2}\\
  &\leq 2 \sqrt{1+\norm{a}^2+\norm{b}^2}
  \leq 2 (1+\norm{a}+\norm{b}). 
\end{align*}
\end{lemma}
Using these inequalities and the 
Poincar\'e inequality we get the following.

\begin{lemma}
The area functional $\AF$ is a priori defined on the function space 
$W^{1,1}(\Omega)$. We can furthermore restrict our attention to 
those functions $u \in W^{1,1}(\Omega)$ with $\int_\Omega u \,\dvol_g=0$,
because $\AF(u)=\AF(u+\const)$.
\end{lemma}

The Banach space $W^{1,1}(\Omega)$ is not weakly sequentially compact. But
$\AF$ can be naturally extended to the space of functions of bounded variation 
$\BV(\Omega)$, the bidual of $W^{1,1}(\Omega)$.
Sets of functions uniformly bounded in the
$\BV$--norm are relatively compact in $L^1(\Omega)$.
It is this property that makes $\BV(\Omega)$ the suitable space for 
investigating variational problems corresponding to (area-like) functionals
with linear growth. For more information on $\BV$--functions see
\cite{Giusti1}, for example.

Extending $\AF$ to $\BV(\Omega)$, minimizing sequences have weak limits in
$\BV(\Omega)$ as usual, compare \cite{GMS} for this approach. 
The drawback is that
$\BV$--functions can and do in general have jumps. We have to use
other arguments to exclude this behaviour. The trick is to prove 
convergence in $W^{1,\infty}$ of a limiting sequence.
For this we do not necessarily need the space $\BV(\Omega)$.

We now collect some of the main features of the area functional $\AF$ and its
integrand $f$:

\begin{lemma}
\label{62}
\begin{enumerate}
\item $f:T_x\Omega \to \bbR$ is Lipschitz continuous for every $x \in \Omega$
  with Lipschitz constant $2k(x)$.
\item $f:T_x\Omega \to \bbR$ is strictly convex for every $x \in \Omega$,
  i.~e., for all $t \in (0,1)$ and $X,Y \in T_x\Omega$, $X \neq Y$:
  $f(tX+(1-t)Y) < tf(X) + (1-t) f(Y)$.
\item A minimum of $\AF$ is unique (up to constants).
\item The functional $\AF$ is sequentially lower semicontinuous with respect
  to weak convergence in $W^{1,p}_{\loc}(\Omega)$, $1 \leq p < \infty$.
\end{enumerate}
\end{lemma}

\begin{proof}
(1) and (2) are easy computations. (3) is a direct consequence of the
    convexity property (2). For (4) see \cite[p.~22, Theorem 2.5]{Giaquinta1}. 
\end{proof}


\subsection{The Euler-Lagrange equation}
\label{56}


The next step is to introduce the Euler-Lagrange equation corresponding to our
variational problem $\AF(u)=\min$ 
and to study its main properties.

\begin{definition}
Let $u:\Omega \to \bbR$ be fixed.
Corresponding to the integrand $f:T\Omega \to \bbR$ of the area functional
$\AF$ we define the $1$--form $\df$ and the symmetric $2$--form
$\ddf$ on $\Omega$ by
\begin{align*}
  \df(X) &:=\ddt f({\grad u} + tX)_{|t=0}\\[0.5ex]
  \ddf(X,Y) &:= \tfrac{\partial^2}{\partial t \,\partial s}
                f ({\grad u} + sX + tY)_{|s=t=0}.
\end{align*}
We also write $\df(u;X)$ and $\ddf(u;X,Y)$ to emphasize the dependence on $u$.
\end{definition}

\begin{lemma}
The Euler-Lagrange equation corresponding to the variational problem
$\AF(u)=\min$ is given by
\begin{equation*}
  \Div \df = \Div \df(u;\cdot) = 0.
\end{equation*}
\end{lemma}

\begin{proof}
Suppose $u$ is a smooth minimum of $\AF$.
For $\varphi \in C_0^\infty(\Omega)$ we get
\begin{align*}
  0 &= \ddt \AF(u+t\varphi)_{|t=0} 
     = \int_\Omega \df(u; \grad \varphi) \,\dvol_g
     = -\int_\Omega \varphi \Div(\df)\,\dvol_g.
\end{align*}
As this holds for arbitrary $\varphi \in C_0^\infty(\Omega)$ the claim
follows.
\end{proof}

Unfortunately, this Euler-Lagrange equation is not uniformly elliptic.
Nevertheless, we have estimates for $\ddf$ which
allow existence and regularity theorems for our variational problem.

\begin{lemma}
For simplyfing notation from now on
\begin{align*}
  V &:= k\cdot (W + \grad u),\\
  V^\pm &:= k \cdot (W \pm \grad h + \grad u).
\end{align*}
Furthermore we will write $\scp{\cdot}{\cdot}:=g(\cdot,\cdot)$
for the Riemannian metric on the orbit
space $M^{n-1}$ and $\norm{\cdot}=\norm{\cdot}_g$ for the corresponding norm.
Using these abbreviations, the following holds:
\begin{align*}
  \df(X) &= \frac{k\scp{V^+}{X}}{\sqrt{1+\norm{V^+}^2}}
            + \frac{k\scp{V^-}{X}}{\sqrt{1+\norm{V^-}^2}}\\
  \intertext{and}
  \ddf(X,Y) &= \frac{k^2 \scp{X}{Y}}{\sqrt{1+\norm{V^+}^2}}
               -\frac{k^2\scp{V^+}{X}\scp{V^+}{Y}}
                     {(1+\norm{V^+}^2)^{\frac{3}{2}}}\\
            &\quad + \frac{k^2 \scp{X}{Y}}{\sqrt{1+\norm{V^-}^2}}
               -\frac{k^2\scp{V^-}{X}\scp{V^-}{Y}}
                     {(1+\norm{V^-}^2)^{\frac{3}{2}}}.
\end{align*}
\end{lemma}

\begin{proof}
Straightforward computation.
\end{proof}

A unit normal field on the surface 
$\psi_u(\Omega) = \{\tau(u(x),\sigma(x))\;|\;x \in \Omega\} \subset \hM^n$
is given by
\begin{equation*}
  \nu := \frac{-k \Hor(V)_{|\psi_u}+K_{|\psi_u}}
         {k\sqrt{1+\norm{V}^2}}.
\end{equation*}
Analogously
\begin{equation*}
  \nu^{\pm} := \frac{-k \Hor(V^{\pm})_{|\psi_u}+K_{|\psi_u}}
                    {k\sqrt{1+\norm{V^{\pm}}^2}}
\end{equation*}
defines normal fields on $\psi_{u+h}(\Omega)$ and $\psi_{u-h}(\Omega)$,
the boundary of the set $\hOmega \subset \hM^n$
we want to symmetrize, compare \eqref{12}.

\begin{definition}
\label{20}
We define the ``projections''
$P: T\Omega \to T(\psi_u(\Omega)) \subset T\hM^n$ and
$P^{\pm}: T\Omega \to T(\psi_{u\pm h}(\Omega)) \subset T\hM^n$
by
\begin{align*}
  P(X) &:= \Hor(X)_{|\psi_u}
                      -\left\langle{\Hor(X)_{|\psi_u}},
                          {\nu_{|\psi_u}}\right\rangle_{\hg}
			  \cdot \nu_{|\psi_u}\\
  P^{\pm}(X) &:= \Hor(X)_{|\psi_{u\pm h}}
          -\left\langle{\Hor(X)_{|\psi_{u\pm h}}},
            {\nu^{\pm}_{|\psi_{u\pm h}}}\right\rangle_{\hg}
            \cdot\nu^{\pm}_{|\psi_{u\pm h}}.
\end{align*}
\end{definition}
Using Pythagoras theorem we get
\begin{align*}
  \norm{P(X)}^2_{\hg} 
  = \norm{X}^2 - \tfrac{\scp{V}{X}^2}{1+\norm{V}^2}
  \quad \text{and} \quad
  \norm{P^{\pm}(X)}^2 
  = \norm{X}^2 - \tfrac{\scp{V^{\pm}}{X}^2}{1+\norm{V^{\pm}}^2}.
\end{align*}
Consequently
\begin{equation*}
  \ddf(u; X,X) = k^2 \frac{\norm{P^+(X)}^2}{\sqrt{1+\norm{V^+}^2}}
                +k^2 \frac{\norm{P^-(X)}^2}{\sqrt{1+\norm{V^-}^2}}.
\end{equation*}

\begin{lemma}
For $\ddf$ we have the following estimate:
\begin{equation}
\label{6}
  \mu_1 \cdot \frac{\norm{P(X)}^2_{\hg}}
              {\sqrt{1+\norm{ \grad u}^2}}
  \leq \ddf(u; X,X)
  \leq   \mu_2 \cdot \frac{\norm{P(X)}^2_{\hg}}
              {\sqrt{1+\norm{ \grad u}^2}}
\end{equation}
The constants $\mu_1$, $\mu_2$ depend on the functions $k$ and $h$. More
precisely:
\begin{align*}
  \mu_1 &= \tfrac{k}{2\sqrt{1+k^2\norm{W}^2} 
                   (1+k^2\norm{\grad h}^2)^{\frac{3}{2}}} \\[1ex]
  \mu_2 &= 16k^2\sqrt{1+k^2\norm{W}^2}(1+k^2\norm{\grad h}^2)^\frac{5}{2}.
\end{align*}
\end{lemma}

\begin{proof}
This is a straightforward computation:
It is easy to see
\begin{align*}
  1+\norm{V^\pm}^2 &\leq 4k^2(1+k^2\norm{W}^2)(1+k^2\norm{\grad h}^2)
                             (1+\norm{\grad u}^2), \text{ and}\\
  1+\norm{\grad u}^2 &\leq 4(1+k^2\norm{W}^2)(1+k^2\norm{\grad h}^2)
                            (1+\norm{V^\pm}^2).
\end{align*}
A short calculation then yields
\begin{align*}
  \ddf(X,X) 
  & \geq \frac{k^2\left( \norm{P^+(X)}^2_{\hg}
                        +\norm{P^-(X)}^2_{\hg} \right)}
         {2k\sqrt{(1+k^2\norm{W}^2)(1+k^2\norm{\grad h}^2)
                  (1+\norm{\grad u}^2)}}\\[1ex]
  &\geq \mu_1 \cdot \frac{\norm{P(X)}^2_{\hg}}
              {\sqrt{1+\norm{ \grad u}^2}}
\end{align*}
For the second inequality one computes
\begin{align*}
  \ddf(X,X) 
  &\leq \tfrac{2k^2\sqrt{(1+k^2\norm{W}^2)(1+k^2\norm{\grad h}^2)}}
              {\sqrt{1+\norm{\grad u}^2}}
        \left(\norm{P^+(X)}^2_{\hg} +\norm{P^-(X)}^2_{\hg} \right)\\[0.5ex]
  &\leq   \mu_2 \cdot \frac{\norm{P(X)}^2_{\hg}}
              {\sqrt{1+\norm{ \grad u}^2}}.
\end{align*}
\end{proof}

Summarizing, instead of uniform ellipticity we only have inequality \eqref{6}.
A priori gradient estimates for $C^2$--solutions
of partial differential equations
which fulfill such an inequality have already been obtained by
Ladyzhenskaya and Ural'tseva \cite{LU1}.
Giaquinta, Modica, and Soucek then used these a priori estimates in order
to obtain existence and regularity results for Dirichlet boundary value
problems corresponding to functionals with linear growth \cite{GMS}.
We will apply these ideas to our minimizing problem $\AF(u)=\min$.


\subsection{Approximating minimizing problems}
\label{57}


One of the main properties of $\BV$--functions is that they may have jumps.
In our situation we would not like a minimum of the area functional $\AF$
to have jumps, because that would mean our symmetrization procedure could tear
sets $\hOmega \subset \hM^n$ apart.
In order to remedy this problem we need gradient estimates for minimzers
of $\AF$. 
It turns out that this can be done without extending $\AF$ to 
$\BV(\Omega)$, we can simply use Sobolev spaces.
For this purpose we will consider the approximating functionals 
$\AF_{\varepsilon}$ defined by
\begin{align*}
  \AF_{\varepsilon}(u) = \int_{\Omega} f_{\varepsilon}(\grad u) \,\dvol_g
  &:=  \int_{\Omega} f(\grad u) + \varepsilon \norm{\grad u}^2 \,\dvol_g
\end{align*}
with $\varepsilon > 0$ and
\begin{enumerate}
\item $u \in W^{1,2}(\Omega) \subset W^{1,1}(\Omega)$, because
  $\vol(\Omega) < \infty$,
  \smallskip
\item $\int_{\Omega} u \,\dvol_g = 0$ as a normalization, because
  $\AF_{\varepsilon}(u) = \AF_{\varepsilon}(u + \const)$.
\end{enumerate}

\begin{definition}
Let $u:\Omega \to \bbR$ be fixed. Corresponding to the integrand
$f_{\varepsilon}:T\Omega \to \bbR$  of the functional $\AF_{\varepsilon}$
we define the $1$--form $\df_{\varepsilon}$ and the symmetric $2$--form 
$\ddf_{\varepsilon}$ on $\Omega$ by
\begin{align*}
  \df_{\varepsilon}(X) &:=\ddt f_{\varepsilon}
                               ({\grad u} + tX)_{|t=0}\\[0.5ex]
  \ddf_{\varepsilon}(X,Y) &:= \tfrac{\partial^2}{\partial t \,\partial s}
                f_{\varepsilon} ({\grad u} + sX + tY)_{|s=t=0}.
\end{align*}
We also write $\df_{\varepsilon}(u;X)$ and 
$\ddf_{\varepsilon}(u;X,Y)$ to emphasize the dependence on  $u$.
\end{definition}

$\df_{\varepsilon}$ and $\ddf_{\varepsilon}$ can easily be computed as
\begin{align*}
  \df_{\varepsilon}(u; X)   &= \df(u; X) 
                               + 2\varepsilon \scp{\grad u}{X}\\[0.5ex]
  \ddf_{\varepsilon}(u; X,Y) &= \ddf(u; X,Y) + 2\varepsilon \scp{X}{Y}.
\end{align*}
Furthermore, the functional $\AF_{\varepsilon}$ corresponds to 
the uniformly elliptic Euler-Lagrange equation
\begin{equation}
\label{30}
\Div \df_{\varepsilon} = \Div \df_{\varepsilon}(u;\cdot) = 0.
\end{equation}
Therefore the standard theory of elliptic partial differential equations
yields

\begin{lemma}
For every $\varepsilon > 0$ the variational problem 
$\AF_{\varepsilon}(u)= \min$ has a unique solution $u_{\varepsilon}$
such that 
\begin{equation*}
  u_{\varepsilon} \in W^{1,2}(\Omega) \cap C^{\infty}(\Omega)
  \text{ and }
  \int_\Omega u_{\varepsilon} \,\dvol_g = 0.
\end{equation*}
\end{lemma}

For any $\varphi \in W^{1,2}(\Omega)$, with 
$\int_{\Omega} \varphi \,\dvol_g = 0$ we have
\begin{equation*}
\begin{split}
  \AF(u_{\varepsilon}) + \varepsilon \int_{\Omega} 
  \norm{\grad u_{\varepsilon}}^2 \,\dvol_g
  \leq \AF_{\varepsilon}(\varphi)
  &\leq \AF_1(\varphi) = \const < \infty.
\end{split}
\end{equation*}
Therefore we have $\varepsilon \int_{\Omega} 
\norm{\grad u_\varepsilon}^2 \,\dvol_g \leq \const < \infty$
for $0< \varepsilon \leq 1$ and by Lemma \ref{10} also
$\int_{\Omega} \norm{\grad u_{\varepsilon}}\,\dvol_g \leq 
 \AF(u_{\varepsilon}) \leq \const < \infty$
for $0<\varepsilon \leq 1$.

Consequently, applying Poincar\'e's inequality  yields for
all $0<\varepsilon \leq 1$
\begin{equation}
\begin{split}
\label{49}
  \varepsilon \int_{\Omega} 
    \abs{u_\varepsilon}^2 \,\dvol_g \leq \const < \infty
  \quad \text{and} \quad
  \int_{\Omega} \abs{u_{\varepsilon}}\,\dvol_g
   \leq \const < \infty.
\end{split}
\end{equation}
Suppose now that for every $U \subset\subset \Omega$ and 
$0<\varepsilon \leq 1$ we have estimates
\begin{equation}
\label{11}
\begin{split}
  \sup_U \abs{u_{\varepsilon}} & \leq C_U < \infty, \text{ and}\\
  \sup_U \norm{\grad u_{\varepsilon}} & \leq C_U < \infty,
\end{split}
\end{equation}
where $C_U$ is a constant depending on $U$ but independent of $\varepsilon$.
Then the existence of a locally uniformly convergent subsequence 
$u_{\varepsilon_i} \to u_0$ follows from the Arzela-Ascoli theorem.
$u_0$ is a locally Lipschitz continuous function.

Since $W^{1,2}(U)$ is weakly sequentially compact for 
$U \subset \subset \Omega$, we can assume $u_{\varepsilon_i} \rightharpoonup
u_0$ weakly in $W^{1,2}_{\loc}(\Omega)$. 
Weakly lower semicontinuity of $\AF$ in $W^{1,2}_{\loc}(\Omega)$ implies 
$\AF(u_0) \leq \liminf_{i \to \infty} \AF(u_{\varepsilon_i})$.
Hence $\AF_{\varepsilon_i}(u_{\varepsilon_i})
\leq\AF_{\varepsilon_i}(\varphi)$ 
implies 
\begin{equation}
\label{8}
 \AF(u_0) \leq \AF(\varphi) \quad 
 \text{for all }  \varphi \in W^{1,2}(\Omega).
\end{equation}
We consider now the minimizing problem 
\begin{equation}
\label{9}
 \AF(u) \to \min, \quad u \in W^{1,1}(\Omega),\quad\int_\Omega u\,\dvol_g=0.
\end{equation}
$\partial \Omega$ is smooth and therefore any $u \in W^{1,1}(\Omega)$ can be
approximated in the $W^{1,1}(\Omega)$--norm by a sequence 
$\varphi_j \in C^{\infty}(\overline{\Omega}) \subset W^{1,2}(\Omega)$.
Since by Lemma \ref{62} $f$ is Lipschitz continuous we have
\begin{align*}
  \abs{\AF(u)-\AF(\varphi_j)} 
   &\leq  \int_{\Omega} 2k\, \|\grad u-\grad \varphi_j\| \,\dvol_g \to 0
\end{align*}
for $j \to \infty$. In other words $\AF(u)=\lim_{j \to \infty} \AF(\varphi_j)$
for $\varphi_j \to u$ in $W^{1,1}(\Omega)$.
Henceforth we know that for the minimizing problem
\eqref{9} there exists a minimizing sequence 
$\varphi_j \in C^{\infty}(\overline{\Omega}) \subset W^{1,2}(\Omega)$.
Inserting this sequence into inequality (\ref{8}) yields
$\AF(u_0) \leq \AF(u)$ for all $u \in W^{1,1}(\Omega)$.

Now the locally uniform convergence of $u_{\varepsilon_i} \to u_0$ and 
$\int_{\Omega} u_{\varepsilon_i} \,\dvol_g = 0$ 
imply $\int_{\Omega} u_0 \,\dvol_g = 0$.
Furthermore Lemma \ref{10} yields
\begin{equation*}
 \int_{\Omega} \|\grad u_0\| \,\dvol_g \leq \AF(u_0) \leq \const < \infty.
\end{equation*}
By Poincar\'e's inequality 
$\int_{\Omega} |u_0| \leq \const < \infty$.
Summarizing we have

\begin{proposition}
\label{39}
Suppose estimates \eqref{11} hold for any $U \subset \subset \Omega$ and
$0 < \varepsilon \leq 1$. Then the minimizing problem $\AF(u) \to \min$, 
$u \in W^{1,1}(\Omega)$, $\int_{\Omega} u \,\dvol_g = 0$ has a locally
Lipschitz continuous solution
\[ u_0 \in W^{1,1}(\Omega) \cap C^{0,1}_{\loc}(\Omega).\]
As the integrand of $\AF$ is strictly convex and independent of the value of
$u$, this is even the unique solution of the minimizing problem. 
\end{proposition}

\begin{remark}
Applying standard regularity theory for elliptic partial 
differential equations of second order we get
\[ u_0 \in C^{\infty}(\Omega)\]
For a short overview of regularity theory see \cite[Appendix C]{Giusti1}.
\end{remark}

\begin{remark}
What remains to be done is to show that the estimates \eqref{11} hold for any
$U \subset \subset \Omega$ and $0 < \varepsilon \leq 1$. 
Proving $\sup_U \abs{u_{\varepsilon}} \leq C_U < \infty$
can be accomplished using a hair cutting argument, see 
\cite[Thm.~14.10, p.~167]{Giusti1} or \cite[Section 3.5]{John1}.
As this is quite standard, we will omit the proof. 
It turns out that the constant $C_U$ in this inequality only depends on
$\dist(\partial U,\Omega)$, $\int_{\Omega} \abs{u_{\varepsilon}}\,\dvol_g$ and 
$\varepsilon \int_{\Omega} \abs{u_{\varepsilon}}^2\,\dvol_g$.
By \eqref{49} we already know that 
$\int_{\Omega} \abs{u_{\varepsilon}} \dvol_g$ and 
$\varepsilon \int_{\Omega} \abs{u_{\varepsilon}}^2 \dvol_g$ can easily be
estimated by constants independent of $\varepsilon$.
\smallskip

The difficult part is to show the estimate
$\sup_U \norm{\grad u_{\varepsilon}} \leq C_U < \infty$.
This is the objective of Section \ref{37}.
\end{remark}


\subsection{Alternative symmetrization procedures}


The initial idea for our symmetrization procedure was to deform a given set
$\hOmega \subset \hM^n$ along the integral lines of the Killing field
$K$ corresponding to a transvection $\tau$, compare Section \ref{4}.
A question that naturally arises in this context is the following:
Besides Killing fields, which other vector fields $X$ could be used to
establish symmetrization procedures?

As we want the symmetrization procedure to be volume preserving, an immediate
consequence is to restrict attention to vector fields $X$ that are 
divergence free. 
Considering an arbitrary vector field of this kind, the resulting
symmetrization procedure will typically have some major (analytical)
disadvantages. To illustrate this just take the vector field
$\grad \beta_{\gamma}$, where $\beta_{\gamma}$ is a Busemann function.
Rescaling, we easily obtain a divergence free vector field
$X=(\varphi \circ \beta_{\gamma}) \cdot \grad \beta_{\gamma}$, where
$\varphi \in C^{\infty}(\bbR)$ is an appropriate function.
In this setting the area functional that corresponds to the area functional in
Lemma \ref{50} takes the form
\begin{equation*}
  \area(u) = 
  \int_\Omega \sqrt{u^2 + \frac{1}{m^2}
  \sum_{i=1}^{n-1} \frac{du(X_i)^2}{u^\frac{2\lambda_i}{m}}}\,\dvol,
\end{equation*}
where $\Omega$ is an open subset of a horosphere and
the $X_i$ are an orthonormal frame of eigenvectors corresponding to 
the eigenvalues $\lambda_i \geq 0$ of the second fundamental form for the
horosphere. 

The first problem here is that the area functional depends on the values of
$u$ directly. This causes severe difficulties aready for the existence 
problem for weak minimizers.
Another problem emerges from the fact that the integrand is not convex in
$(u,du)$, making this area functional more or less useless for our purposes.


\section{The gradient estimate}
\label{37}


The core of this section is to establish the estimate
\begin{equation*}
  \sup_U \norm{\grad u_{\varepsilon}} \leq C_U < \infty, 
  \quad U \subset\subset \Omega, \quad 0 < \varepsilon \leq 1,
\end{equation*}
for the gradient of the minimizer $u_{\varepsilon}$ of $\AF_{\varepsilon}$.
For this purpose we will introduce a differential equation for 
$\grad u_{\varepsilon}$, which will subsequently be
used to estimate this gradient.


\subsection{Differential equation for the gradient}


A differential equation for $\grad u_{\varepsilon}$ is obtained by 
differentiating the Euler-Lagrange equation for the minimizer 
$u_{\varepsilon}$ of $\AF_{\varepsilon}$ and exchanging the order of
differentiation.

\begin{definition}
From now on, we will use the abbreviation
\begin{equation*}
  p_{\varepsilon} := \norm{\grad u_{\varepsilon}}^2.
\end{equation*}
\end{definition}

Observe that we have the following identity
\begin{align}
\label{14}
  \grad p_{\varepsilon} = 2 \nab_{\grad u_{\varepsilon}}\grad u_{\varepsilon}
  = 2 \Hess u_{\varepsilon} (\grad u_{\varepsilon}).
\end{align}

\begin{lemma}
\label{13}
Let $X_1,\dotsc, X_{n-1}$ be an arbitrary orthonormal frame on $\Omega$.
For $X \in T\Omega$ and the $1$--form $\df_{\varepsilon}$ (or any other
$1$--form instead) we have
\begin{equation*}
     d_X(\Div  \df_{\varepsilon}) 
  = \Div(\nab_X \df_{\varepsilon})-
        (\nabla \df_{\varepsilon})(\nab_{X_j}X,X_j)
     -\df_{\varepsilon} \cdot \Ric(X).
\end{equation*}
\end{lemma}

\begin{proof}
Here and in the following computations we use the convention to sum over
repeated indices.
\begin{align*}
      d_X(\Div \df_{\varepsilon}) 
  &= \nab_X((\nab \df_{\varepsilon})(X_j,X_j))
   = (\nab^2_{X,X_j} \df_{\varepsilon})(X_j)\\
  &= (\nab^2_{X_j,X} \df_{\varepsilon})(X_j)
     -\df_{\varepsilon}(R(X,X_j)X_j)\\
  &= \Div(\nab_X \df_{\varepsilon})
     -(\nab \df_{\varepsilon})(\nab_{X_j} X,X_j)
      - \df_{\varepsilon}\cdot \Ric(X).
\end{align*}
\end{proof}

From now on in this chapter, we will always consider $\ddf_{\varepsilon}$
and $\df_{\varepsilon}$ corresponding to the unique minimizer
$u_{\varepsilon}$ of $\AF_{\varepsilon}$, that is
$\ddf_{\varepsilon} = \ddf_{\varepsilon}(u_{\varepsilon};\cdot,\cdot)$,
and $\df_{\varepsilon} = \df_{\varepsilon}(u_{\varepsilon};\cdot)$.

\begin{definition}
For $\varepsilon > 0$ the $2$--form
$\ddb_{\varepsilon}$ on $\Omega$ is defined by
\begin{equation*}
  \ddb_{\varepsilon}(X,Y):=\nab_X \df_{\varepsilon} \cdot Y 
  - \ddf_{\varepsilon}(\nab_X \grad u_{\varepsilon}, Y).
\end{equation*}
\end{definition}

\begin{lemma}
Let $X_1,\dotsc, X_{n-1}$ be an orthonormal frame on $\Omega$.
Then the following differential equation holds for $\grad u_{\varepsilon}$:
\begin{equation}
\label{16}
  0 = \frac{1}{2} \Div\left
        (\ddf_{\varepsilon}(\grad p_{\varepsilon}, \cdot) \right) 
      - \ddf_{\varepsilon}(\nab_{X_j}\grad u_{\varepsilon},
                           \nab_{X_j} \grad u_{\varepsilon})
      + B,
\end{equation}
where  
\begin{equation}
\label{15}
  B := \Div\left(\ddb_\varepsilon(\grad u_{\varepsilon}, \cdot)\right) 
        - \ddb_\varepsilon(\nab_{X_j} \grad u_{\varepsilon}, X_j)
        - \df_{\varepsilon} \cdot \Ric(\grad u_{\varepsilon}).
\end{equation}
\end{lemma}

\begin{proof}
Differentiating the Euler-Lagrange equation 
$\Div \df_{\varepsilon} = 0$ and applying Lemma \ref{13} 
with $X=\grad u_{\varepsilon}$ yields
\begin{align*}
 0 &= d_{\grad u_{\varepsilon}}(\Div \df_{\varepsilon})\\
   &= \Div \left( \ddf_{\varepsilon}(
          \nab_{\grad u_{\varepsilon}}\grad u_{\varepsilon}, \cdot)
        + \ddb_{\varepsilon}(\grad u_{\varepsilon},\cdot)\right)\\
   &  \quad -(\nab\df_{\varepsilon})(\nab_{X_j}\grad u_{\varepsilon},
       X_j)-\df_{\varepsilon}\cdot\Ric(\grad u_{\varepsilon}).
\end{align*}
Now
\begin{align*}
 (\nab \df_{\varepsilon})&(\nab_{X_j}\grad u_{\varepsilon},X_j)\\
 &\begin{aligned}
   &=\ddf_{\varepsilon}(\nab_{X_j}\grad u_{\varepsilon},
                      \nab_{X_j}\grad u_{\varepsilon})
   +\ddb_{\varepsilon}(\nab_{X_j}\grad u_{\varepsilon}, X_j).
 \end{aligned}
\end{align*}
Combining these identities with \eqref{14} gives the 
differential equation.
\end{proof}

\begin{lemma}
For the quantity $B$ defined in \eqref{15} we have the estimate
\begin{equation}
\label{24}
  B \leq \abs{B} \leq \mu_3 \left(\sum_{i=1}^{n-1} 
    \norm{P(\nab_{X_i}\grad u_{\varepsilon})} \right)
    + \mu_4 \sqrt{1+p_{\varepsilon}}
\end{equation}
with $\mu_3=\mu_3(x)$, $\mu_4=\mu_4(x)$.
$X_1,\dotsc, X_{n-1}$ denotes an arbitrary orthonormal frame on $\Omega$.
\end{lemma}

\begin{proof}
This is a long but straightforward calculation.
\end{proof}


\subsection{Preliminary estimates}


For $U \subset \subset \Omega$ we want to estimate
$\sup_U \norm{\grad u_{\varepsilon}}$. Therefore we define

\begin{definition}
For $p_{\varepsilon} = \norm{\grad u_\varepsilon}^2$ 
\begin{equation*}
   q_{\varepsilon} := \log(1+p_{\varepsilon}).
\end{equation*}
\end{definition}

Observe that we have the identities
\begin{align*}
  \grad q_{\varepsilon} = \frac{\grad p_{\varepsilon}}{1+p_{\varepsilon}}
  \quad \text{and} \quad
  P(\grad q_{\varepsilon}) = \frac{P(\grad p_{\varepsilon})}
                                   {1+p_{\varepsilon}}.
\end{align*}

\begin{definition}
We define the sets
\begin{align*}
  \Omega_{\lambda} &:= \{ x \in\Omega \;|\; q_{\varepsilon}(x) > \lambda\} \\
  \Omega_{\lambda,\rho} &:= \Omega_{\lambda} \cap B_\rho(x_0),\\
\end{align*}
where $x_0 \in \Omega$ is an arbitrary but fixed point and $\rho \leq R_0$ with
$R_0$ such that $B_{R_0}(x_0) \subset \subset \Omega$.
Furthermore
\begin{align*}
  S_\lambda &:= 
  \{ (x,u_{\varepsilon}(x)) \in \Omega \times \bbR 
                            \;|\; x \in \Omega_{\lambda} \},\\
  S_{\lambda,\rho} &:= (\Omega_{\lambda,\rho} \times \bbR) \cap S_{\lambda}.
\end{align*}
\end{definition}

The idea for estimating $\sup \norm{\grad u_{\varepsilon}}$ works as follows:
Consider  
\begin{align*}
  \beta(\lambda,\rho) &:= 
     \int_{S_{\lambda,\rho}} (q_{\varepsilon}-\lambda)^2 \,\dHm_{n-1}
     +\varepsilon \int_{\Omega_{\lambda,\rho}} (1+p_{\varepsilon})
       (q_{\varepsilon}-\lambda)^2 \,\dvol_g\\
   &\phantom{:}= \int_{\Omega_{\lambda,\rho}} (q_{\varepsilon}-\lambda)^2
      \cdot \sqrt{1+p_{\varepsilon}} \,\dvol_g
     +\varepsilon \int_{\Omega_{\lambda,\rho}} (1+p_{\varepsilon})
       (q_{\varepsilon}-\lambda)^2 \,\dvol_g
\end{align*}
where $\Hm_{n-1}$ denotes the $(n-1)$--dimensional Hausdorff measure
corresponding naturally to $S_{\lambda,\rho} \subset \Omega \times \bbR$.
Using the differential equation \eqref{16} we will derive 
some estimates for the terms involved in the definition of $\beta$.
These will be applied to show that there exist 
$0 < \rho_0,\lambda_0 < \infty$ such that
$\beta(\lambda_0,\rho_0) = 0$.
This is just equivalent to
\begin{equation*}
 \esssup_{x \in B_{\rho_0}(x_0)} q_{\varepsilon}(x) \leq \lambda_0 < \infty,
\end{equation*}
which is nothing but the desired gradient estimate.\\
We will now establish the basic estimates involved into these computations.

\begin{lemma}
\label{17}
For $p_{\varepsilon} = \norm{\grad u_{\varepsilon}}^2$ we have
\begin{equation*}
 \norm{\grad p_{\varepsilon}}^2 \leq 4 \norm{\grad u_{\varepsilon}}^2\cdot
     \sum_{i=1}^{n-1}\norm{\nab_{X_i}\grad u_{\varepsilon}}^2,
\end{equation*}
where $X_1,\dotsc,X_{n-1}$ is an arbitrary orthonormal frame on $\Omega$.
\end{lemma}

\begin{proof}
This is just an application of the Cauchy-Schwarz inequality.
\end{proof}

\begin{lemma}
\label{18}
Let $X_1,\dotsc,X_{n-1}$ be an orthonormal frame 
of $\Omega$, $P$ the projection as in Definition \ref{20}. Then
\begin{equation}
\label{19}
  \norm{P(\grad p_\varepsilon)}_{\hg}^2 
  \leq  4 \norm{\grad u_{\varepsilon}}^2
        \sum_{i=1}^{n-1}\norm{P(\nab_{X_i} \grad u_{\varepsilon})}^2_{\hg}.
\end{equation}
\end{lemma}

\begin{proof}
It is easy to show that the expressions in \eqref{19} are independent of the
choice of the orthonormal frame.
Therefore we can choose an orthonormal
frame $\{X_i\}$ such that 
$X_1=\frac{\grad u_{\varepsilon}}{\norm{\grad u_{\varepsilon}}}$.
We abbreviate $V_{\varepsilon}:=k \cdot (W+\grad u_{\varepsilon})$.
The Cauchy-Schwarz inequality now implies
\begin{align*}
  &\;
    \Big( 4 p_{\varepsilon}  
    \scp{\nab_{X_1} \grad u_{\varepsilon}}{V_{\varepsilon}}^2
    - \scp{\grad p_{\varepsilon}}{V_{\varepsilon}}^2 \Big)
            + 4 p_{\varepsilon} \sum_{i=2}^{n-1} 
     \scp{ \nab_{X_i} \grad u_{\varepsilon}}{V_{\varepsilon}}^2\\
  \leq&\;
    (1+\norm{V_{\varepsilon}}^2) 
       \Big(4p_{\varepsilon} \norm{\nab_{X_1}\grad u_{\varepsilon}}^2
      -\norm{\grad p_{\varepsilon}}^2
    +4 p_{\varepsilon} 
        \sum_{i=2}^{n-1}\norm{\nab_{X_i} \grad u_{\varepsilon}}^2
       \Big).
\end{align*}
An easy computation yields the claim.
\end{proof}

\begin{lemma}
\label{21}
Let $\zeta \in C_0^\infty(\Omega)$, $\lambda \geq 0$. Then
\begin{equation}
\begin{split}
\label{22}
  &\int_{S_\lambda} \norm{P(\grad q_{\varepsilon})}^2\zeta^2 \,\dHm_{n-1}
  + \varepsilon \int_{\Omega_\lambda} 
  (1+p_{\varepsilon})\norm{\grad q_{\varepsilon}}^2\zeta^2\,\dvol_g\\
  \leq & \; C \cdot \Big\{ 
     \int_{S_{\lambda}} (q_{\varepsilon}-\lambda)^2
     \norm{P(\grad \zeta)}^2
          + (q_{\varepsilon}-\lambda) \zeta^2 \,\dHm_{n-1} \\
  & \qquad + \varepsilon \int_{\Omega_\lambda} 
    (1+p_{\varepsilon})(q_{\varepsilon}-\lambda)^2
               \norm{\grad \zeta}^2 \,\dvol_g \Big\}.
\end{split}
\end{equation}
$C$ depends on $\inf_{\supp \zeta} \mu_1$,
$\sup_{\supp \zeta} \mu_2$,
$\sup_{\supp \zeta} \mu_3$, and
$\sup_{\supp \zeta} \mu_4$.
\end{lemma}

\begin{proof}
Choose the test function
\begin{equation*}
  \varphi(x) := \zeta^2(x) \cdot \max\{q_{\varepsilon}(x)-\lambda,0\}, 
\end{equation*}
with $\zeta \in C_0^\infty(\Omega)$.
Multiplying the differential equation \eqref{16} with $\varphi$, integrating
over $\Omega$ and applying the divergence theorem yields
\begin{equation}
\begin{split}
\label{23}
  \int_{\Omega_{\lambda}} &
  \underbrace{
    \tfrac{1}{2} \ddf_{\varepsilon}(\grad p_{\varepsilon},
                                    \grad q_{\varepsilon}) 
    \cdot \zeta^2}_{(1.)}\\
  & + \underbrace{
    \ddf_{\varepsilon}(\nab_{X_i} \grad u_{\varepsilon},
                       \nab_{X_i} \grad u_{\varepsilon})
    \cdot \zeta^2\cdot (q_{\varepsilon}-\lambda)}_{(2.)} \dvol_g \\
  = \int_{\Omega_{\lambda}} &
     \underbrace{
      -\ddf_{\varepsilon}(\grad p_{\varepsilon},
                          \grad \zeta) 
       \cdot \zeta \cdot (q_{\varepsilon}-\lambda)}_{(3.)}
     + \underbrace{
     B \cdot \zeta^2\cdot (q_{\varepsilon}-\lambda)}_{(4.)} \dvol_g
\end{split}
\end{equation}
We now estimate the expressions (1.) -- (4.) separately.\\
For (1.) we apply \eqref{6} to obtain
\begin{align*}
    \tfrac{1}{2} \ddf_{\varepsilon}&(\grad p_{\varepsilon},
                                     \grad q_{\varepsilon}) 
                  \cdot \zeta^2 \\
   &\begin{aligned}
    &= \tfrac{1}{2} \left(\ddf(\grad q_{\varepsilon},
	    	               \grad q_{\varepsilon})
       + 2\varepsilon \scp{\grad q_{\varepsilon}}{\grad q_{\varepsilon}}
       \right) (1+p_{\varepsilon}) \zeta^2\\
    &\geq
       \tfrac{1}{2} \mu_1 \norm{P(\grad q_{\varepsilon})}^2 
       \sqrt{1+p_{\varepsilon}} \cdot \zeta^2
       + \varepsilon \,\norm{\grad q_{\varepsilon}}^2
                  (1+p_{\varepsilon})\zeta^2.
   \end{aligned}
\end{align*}
For (2.) we use \eqref{6} again to compute
\begin{align*}
  \ddf_{\varepsilon}(\nab_{X_i} \grad u_{\varepsilon},
                     \nab_{X_i} \grad u_{\varepsilon})
  &\geq
      \mu_1 \frac{\sum_i \norm{P(\nab_{X_i}\grad u_{\varepsilon})}^2}
                       {\sqrt{1+p_{\varepsilon}}}. 
\end{align*}
As $f_{|T_x\Omega}$ is strictly convex for every $x \in \Omega$, we know that
$\ddf$ is positive definit. Therefore the Cauchy-Schwarz
inequality and \eqref{6} applied to (3.) gives
\begin{align*}
  -\ddf_\varepsilon&(\grad p_{\varepsilon},
                        \grad \zeta)  \\[0.5ex]
  &\begin{aligned}
    &\leq
      \mu_2 
        \frac{\norm{P(\grad p_{\varepsilon})}\norm{P(\grad \zeta)}}
             {\sqrt{1+p_{\varepsilon}}} 
         + 2\varepsilon \norm{\grad p_{\varepsilon}}
                        \norm{\grad \zeta}.  
  \end{aligned}
\end{align*}
(4.) can obviously be estimated by \eqref{24}.\\
We define the constants 
\begin{equation*}
\begin{aligned}
   C_1 &:= \inf_{\supp \zeta} \mu_1 \quad \text{and} \quad 
   C_i &:= \sup_{\supp \zeta} \mu_i, \quad i=2,3,4.
\end{aligned}
\end{equation*}
Now we insert the estimates for (1.) -- (4.) into (\ref{23}) to calculate
\begin{align*}
  &
  \begin{aligned}
    &   \frac{C_1}{2} \int_{S_{\lambda}} 
        \norm{P(\grad q_{\varepsilon})}^2 \zeta^2 \,\dHm_{n-1}\\
    & + C_1 \int_{S_{\lambda}}
        \frac{\sum_i \norm{P(\nabla_{X_i} \grad u_\varepsilon)}^2}
            {1+p_{\varepsilon}}
        \zeta^2\,(q_{\varepsilon}-\lambda) \,\dHm_{n-1}\\
    & + \varepsilon \int_{\Omega_{\lambda}}
        \norm{\grad q_{\varepsilon}}^2(1+p_{\varepsilon})
        \zeta^2 \,\dvol_g\\ 
  \end{aligned}\\[1ex]
  \displaybreak[0]
  \leq &
  \begin{aligned}[t]
    &   C_2 \int_{S_{\lambda}} 
        \frac{\norm{P(\grad p_{\varepsilon})}\norm{P(\grad \zeta)}}
             {1+p_{\varepsilon}}
        \zeta \, (q_{\varepsilon}-\lambda) \,\dHm_{n-1}\\
    & + C_3 \int_{S_{\lambda}} 
        \frac{\sum_i \norm{P(\nabla_{X_i}\grad u_\varepsilon)}}
             {\sqrt{1+p_{\varepsilon}}}
        (q_{\varepsilon}-\lambda) \zeta^2\,\dHm_{n-1}\\[0.5ex]
    & + C_4 \int_{S_{\lambda}}
            (q_{\varepsilon}-\lambda)\zeta^2 \,\dHm_{n-1}\\[0.5ex]
    & + 2 \varepsilon \int_{\Omega_{\lambda}}
        \norm{\grad p_{\varepsilon}} \norm{\grad \zeta}\,
        \zeta\,(q_{\varepsilon}-\lambda)\,\dvol_g
  \end{aligned}
  \displaybreak[0]
\intertext{Using H\"older's inequality and the Cauchy inequality 
$ab \leq \varepsilon_l \,a^2 + \tfrac{1}{4\varepsilon_l}\,b^2$,
$l=1,2,3$, we continue}
 \leq &
  \begin{aligned}[t]
    &   C_2 \varepsilon_1 \int_{S_{\lambda}} 
                          \norm{P(\grad q_{\varepsilon})}^2\,
                          \zeta^2 \,\dHm_{n-1}\\
    & + \frac{C_2}{4\varepsilon_1}
        \int_{S_{\lambda}}\norm{P(\grad \zeta)}^2\,
        (q_{\varepsilon}-\lambda)^2 
        \,\dHm_{n-1}\\
    & + C_3 \varepsilon_2 
        \int_{S_{\lambda}} 
        \frac{\sum_i \norm{P(\nabla_{X_i}\grad u_\varepsilon)}^2}
             {1+p_{\varepsilon}}
        (q_{\varepsilon}-\lambda) \zeta^2 \,\dHm_{n-1}\\
    &  + \left(\frac{C_3}{4 \varepsilon_2}+C_4 \right)
        \int_{S_{\lambda}} (q_{\varepsilon}-\lambda)\zeta^2 \,\dHm_{n-1}\\
    & + 2 \varepsilon \varepsilon_3
        \int_{\Omega_{\lambda}} 
        \norm{\grad q_{\varepsilon}}^2\zeta^2(1+p_{\varepsilon})\,\dvol_g\\
    & + \frac{\varepsilon}{2\varepsilon_3}     
        \int_{\Omega_{\lambda}} \norm{\grad \zeta}^2 
        (q_{\varepsilon}-\lambda)^2(1+p_{\varepsilon})\,\dvol_g.
  \end{aligned}
\end{align*}
Set $C_2\varepsilon_1 =\frac{C_1}{4}$, $C_3 \varepsilon_2=C_1$ and
$\varepsilon_3=\frac{1}{4}$. Performing a short calculation and 
choosing $C$ appropriately,
inequality (\ref{22}) follows.
\end{proof}

\begin{lemma}
\label{25}
Let $\zeta \in C_0^\infty(\Omega)$, $\lambda \geq 0$. Then
\begin{equation}
\begin{split}
\label{26}
  \varepsilon \int_{\Omega_{\lambda}} \sum_{i=1}^{n-1} 
    \norm{\nabla_{X_i} & \grad u_{\varepsilon}}^2
    (q_{\varepsilon}-\lambda)^2 \zeta^2 \dvol_g\\
  &\leq  
  \begin{aligned}[t]
        C \cdot \Big\{
        \varepsilon & \int_{\Omega_{\lambda}} (1+p_{\varepsilon})
        (q_{\varepsilon}-\lambda)^2\norm{\grad \zeta}^2 \dvol_g\\
     +  & \int_{S_\lambda}(q_{\varepsilon}-\lambda)^2
        (\norm{P(\grad \zeta)}^2+\zeta^2) \dHm_{n-1}
  \Big\}.
  \end{aligned}
\end{split}
\end{equation}
$C$ depends on $\inf_{\supp \zeta} \mu_1,\;$
$\sup_{\supp \zeta} \mu_2,\;$
$\sup_{\supp \zeta} \mu_3$, and
$\sup_{\supp \zeta} \mu_4$.
\end{lemma}

\begin{proof}
Similar to the proof of Lemma \ref{21} we now consider the test function
\begin{align*}
 \varphi(x) &:= \zeta^2 \cdot \max\{q_{\varepsilon}(x)-\lambda,0\}^2.
\end{align*}
Multiplying the differential equation \eqref{16} with $\varphi$, integrating
over $\Omega$ and applying the divergence theorem yields
\begin{equation}
\begin{split}
\label{27}
  \int_{\Omega_\lambda} &
  \underbrace{
    \ddf_{\varepsilon}(\grad p_{\varepsilon},
                     \grad q_{\varepsilon}) 
    \cdot (q_{\varepsilon}-\lambda) \cdot \zeta^2}_{(1.)}\\
  & + \underbrace{
    \ddf_{\varepsilon}(\nab_{X_i} \grad u_{\varepsilon},
                     \nab_{X_i} \grad u_{\varepsilon})
    \cdot \zeta^2\cdot (q_{\varepsilon}-\lambda)^2}_{(2.)} \,\dvol_g \\
  =  \int_{\Omega_\lambda} &
     \underbrace{
      -\ddf_{\varepsilon}(\grad p_{\varepsilon},\grad \zeta) 
       \cdot \zeta \cdot 
      (q_{\varepsilon}-\lambda)^2}_{(3.)}
     + \underbrace{
     B \cdot \zeta^2\cdot (q_{\varepsilon}-\lambda)^2}_{(4.)} \,\dvol_g.
\end{split}
\end{equation}
Again we examine the terms (1.) -- (4.) separately.
The first term is simply estimated by
\begin{align*}
    \ddf_{\varepsilon}&(\grad p_{\varepsilon},\grad q_{\varepsilon}) 
      \cdot (q_{\varepsilon}-\lambda) \cdot \zeta^2  \geq 0 
    \quad \text{on }\Omega_{\lambda}.
\end{align*}
For (2.) we get, using \eqref{6},
\begin{align*}
  \ddf_{\varepsilon}&(\nab_{X_i} \grad u_{\varepsilon}, 
                      \nab_{X_i} \grad u_{\varepsilon})\\
  &\begin{aligned}
    & \geq 2 \frac{\mu_1}{2}
             \frac{\sum_i\norm{P(\nab_{X_i}\grad u_{\varepsilon})}^2}
                  {\sqrt{1+p_{\varepsilon}}}
             + 2\varepsilon\sum_i\norm{\nab_{X_i}\grad u_{\varepsilon}}^2 
  \end{aligned}
\end{align*}
The last two terms are estimated independently:
We use Lemma \ref{18} to obtain
\begin{align*}
     \frac{\sum_i\norm{P(\nab_{X_i} \grad u_{\varepsilon})}^2}
      {\sqrt{1+p_{\varepsilon}}}
     &\geq \frac{1}{4}\norm{P(\grad q_{\varepsilon})}^2\sqrt{1+p_{\varepsilon}}
\end{align*}
and by Lemma \ref{17}
\begin{align*}
  \sum_i \norm{\nab_{X_i}\grad u_{\varepsilon}}^2 
  &\geq \frac{1}{4}\norm{\grad q_{\varepsilon}}^2(1+p_{\varepsilon}).
\end{align*}
Summarizing, (2.) can be estimated as
\begin{align*}
 \ddf_{\varepsilon}&(\nab_{X_i} \grad u_{\varepsilon},
                       \nab_{X_i} \grad u_{\varepsilon})\\
 &\geq
  \begin{aligned}[t]
    &   \frac{\mu_1}{2}\cdot\frac{\sum_i 
        \norm{P(\nab_{X_i}\grad u_{\varepsilon})}^2}{1+p_{\varepsilon}}
        \sqrt{1+p_{\varepsilon}}
     + \frac{\mu_1}{8} \norm{P(\grad q_{\varepsilon})}^2
        \sqrt{1+p_{\varepsilon}}\\
    & + \varepsilon \sum_i \norm{\nab_{X_i}\grad u_{\varepsilon}}^2 
     + \frac{\varepsilon}{4}\norm{\grad q_{\varepsilon}}^2
        (1+p_{\varepsilon}).
  \end{aligned}
\end{align*}
For (3.) inequality \eqref{6} yields
\begin{align*}
  -\ddf_{\varepsilon}&(\grad p_{\varepsilon},\grad \zeta) \\
  &\leq 
          \mu_2 \frac{\norm{P(\grad p_{\varepsilon})}\,\norm{P(\grad \zeta)}}
                     {\sqrt{1+p_{\varepsilon}}}
          +2\varepsilon \norm{\grad p_{\varepsilon}}
          \norm{\grad \zeta}.
\end{align*}
(4.) can again be estimated by \eqref{24}.
\smallskip

These estimates can now be inserted into \eqref{27}. Performing then
the same steps as in the proof of Lemma \ref{21},
that is, applying H\"older's inequality and Cauchy's inequality with
$\varepsilon_l$, we finally get the desired inequality \eqref{26}.
\end{proof}

\begin{lemma}
\label{28}
Let $u_{\varepsilon} \in C^2(\overline{\Omega})$ 
be a solution of the differential equation
(\ref{30}) and $\Omega' \subset \Omega$ such that even the convex hull of
$\Omega'$ is contained in $\Omega$.
Then for any $\varphi \in C^\infty(\overline{\Omega'})$ with
$\varphi_{|\partial \Omega'}=0$
\begin{equation}
\label{29}
  \int_{S'} \varphi^2 \,\dHm_{n-1} \leq C 
  \cdot \Hm_{n-1}^{\frac{2}{n-1}}(S')
  \cdot \int_{S'} \norm{P(\grad \varphi)}^2 \,\dHm_{n-1}.
\end{equation}
Here $S'=\{(x,u_\varepsilon(x)) \in \Omega' \times \bbR \}$ and 
$C$ depends on $\osc_{\Omega'} u_{\varepsilon}$, 
as well as on constants involving the values of $k$, $\norm{W}$, 
$\norm{\grad h}$ on $\Omega$. 
\end{lemma}

\begin{proof}
The proof is based on isoperimetric
inequalities and can be accomplished in analogy to \cite[Lemma 2,
p.~697]{LU1}.
See also \cite[Lemma 3.8]{GMS}.
\end{proof}


\subsection{The final estimates}


We have now established all the inequalities necessary to derive the final
estimate for the function  
\begin{align*}
  \beta(\lambda,\rho) &:= 
     \int_{S_{\lambda,\rho}} (q_{\varepsilon}-\lambda)^2 \,\dHm_{n-1}
     +\varepsilon \int_{\Omega_{\lambda,\rho}} (1+p_{\varepsilon})
       (q_{\varepsilon}-\lambda)^2 \,\dvol_g,
\end{align*}
which can then be used to show that there exist $0<\lambda_0,\rho_0<\infty$
such that $\beta(\lambda_0,\rho_0) = 0$. We already mentioned that this is
equivalent to the desired bound for $\sup_U \norm{\grad u_{\varepsilon}}$.

\begin{definition}
Let $x_0 \in \Omega$ be fixed, $R_0 \in \bbR$, 
$B_{R_0}(x_0) \subset \Omega$.
We define $\zeta_{\rho,R}: \Omega \to \bbR$, $0< \rho < R \leq R_0$, to be
a $C^\infty(\Omega)$--test function with 
$\zeta_{\rho,R}(x)=1$ for $x \in B_{\rho}(x_0)$,
$\zeta_{\rho,R}(x)=0$ for $x \in \Omega \backslash B_R(x_0)$,
$\norm{\grad \zeta_{\rho,R}} \leq \frac{\const}{R-\rho}$.
\end{definition}

\begin{lemma}
For all $0< \rho < R \leq R_0$ we have
\begin{eqnarray}
\label{32}
  \beta(\lambda,\rho) \leq
  C \cdot \Hm_{n-1}^{\frac{2}{n-1}}(S_{\lambda,R})
  \left( \frac{1}{(R-\rho)^2} \beta(\lambda,R) 
  + \Hm_{n-1}(S_{\lambda,R})
  \right).
\end{eqnarray}
The constant $C$ depends on the constants appearing in the inequalities of the
previous section.
\end{lemma}

\begin{proof}
Applying Lemma \ref{28} we estimate
\begin{align*}
  \beta(\lambda,\rho) 
  \leq&
    \begin{aligned}
       \int_{S_{\lambda,R}} (q_{\varepsilon}-\lambda)^2 
         \zeta_{\rho,R}^2\,\dHm_{n-1}
       + \varepsilon \int_{\Omega_{\lambda,R}} 
         (1+p_{\varepsilon})(q_{\varepsilon}-\lambda)^2 
         \zeta_{\rho,R}^2\,\dvol_g
    \end{aligned}\\[1ex]
    \displaybreak[0]
  \leq&
    \begin{aligned}[t]
      C \Hm_{n-1}^{\frac{2}{n-1}}(S_{\lambda,R})
      \Bigg(&\int_{S_{\lambda,R}}(q_{\varepsilon}-\lambda)^2
        \norm{P(\grad \zeta_{\rho,R})}^2 \,\dHm_{n-1}\\
      & + \varepsilon \int_{\Omega_{\lambda,\rho}}(1+p_{\varepsilon})
          (q_{\varepsilon}-\lambda)^2
          \norm{\grad \zeta_{\rho,R}}^2 \,\dvol_g\\
      & + \int_{S_{\lambda,R}} \norm{P(\grad q_{\varepsilon})}^2
          \zeta_{\rho,R}^2 \,\dHm_{n-1}\\
      & + \varepsilon \int_{\Omega_{\lambda,R}} 
          \norm{\grad q_{\varepsilon}}^2
          \zeta_{\rho,R}^2 (1+p_{\varepsilon}) \,\dvol_g\\
      & + \varepsilon \int_{\Omega_{\lambda,R}} (q_{\varepsilon}-\lambda)^2
          \zeta_{\rho,R}^2
          \frac{\norm{\grad p_{\varepsilon}}^2}{1+p_{\varepsilon}} \,\dvol_g
      \Bigg).
    \end{aligned}
\end{align*}
We now insert $\norm{P(\grad \zeta_{\rho,R})}^2 \leq
\norm{\grad \zeta_{\rho,R}}^2 \leq \frac{\mathrm{const}}{(R-\rho)^2}$,
Lemma \ref{21} and Lemma \ref{17} to obtain
\begin{align*}
  \beta(\lambda,\rho) 
  \leq
    \begin{aligned}[t]
      C \Hm_{n-1}^{\frac{2}{n-1}}(S_{\lambda,R}) \cdot 
      &\Bigg(
        \frac{\beta(\lambda,R)}{(R-\rho)^2}
        + \underbrace{\int_{S_{\lambda,R}} (q_{\varepsilon}-\lambda) 
          \zeta_{\rho,R}^2 \,\dHm_{n-1}}_{(1.)}  \\
      & + \underbrace{\varepsilon \int_{\Omega_{\lambda,R}} 
            (q_{\varepsilon}-\lambda)^2 \zeta_{\rho,R}^2
            \sum_{i=1}^{n-1}\norm{\nab_{X_i}\grad u_{\varepsilon}}^2 
            \,\dvol_g}_{(2.)}
      \Bigg)
    \end{aligned}
\end{align*}
Using the H\"older inequality (1.) can be estimated as
\begin{align*}
  (1.)
  &\leq \tfrac{1}{2} 
       \int_{S_{\lambda,R}} (q_{\varepsilon}-\lambda)^2\,
  \zeta_{\rho,R}^2 \,\dHm_{n-1}
       + \tfrac{1}{2} \Hm_{n-1}(S_{\lambda,R}). 
\end{align*}
(2.) can be estimated using Lemma \ref{25}. Summarizing
\begin{align*}
  \beta(\lambda,\rho) 
  \leq
    \begin{aligned}[t]
      C \Hm_{n-1}^{\frac{2}{n-1}}(S_{\lambda,R})
      &\cdot \Bigg(
        \frac{\beta(\lambda,R)}{(R-\rho)^2}+ \Hm_{n-1}(S_{\lambda,R})\\
        & + \varepsilon \int_{\Omega_{\lambda,R}}
            (1+p_{\varepsilon})(q_{\varepsilon}-\lambda)^2
            \norm{\grad \zeta_{\rho,R}}^2 \,\dvol_g \\
        & + \int_{S_{\lambda,R}} (q_{\varepsilon}-\lambda)^2
            (\norm{P(\grad \zeta_{\rho,R})}^2+\zeta_{\rho,R}^2)\,\dHm_{n-1}  
      \Bigg).
    \end{aligned}
\end{align*}
Obviously $\zeta_{\rho,R}^2 \leq \frac{\const}{(R-\rho)^2}$ and consequently
(\ref{32}) follows. 
\end{proof}

\begin{remark}
Obviously $(\lambda,\rho) \mapsto \beta(\lambda,\rho)$ and
$(\lambda,\rho) \mapsto \Hm_{n-1}(S_{\lambda,\rho})$ are 
nonnegative functions which are
monotone increasing in $\rho$ and
monotone decreasing in $\lambda$. 
\end{remark}

\begin{lemma}
For $\lambda < \Lambda$ and $R \leq R_0$
\begin{equation}
\label{33}
  \Hm_{n-1}(S_{\Lambda,R})
  \leq \frac{\beta(\lambda,R)}{(\Lambda-\lambda)^2}.
\end{equation}
\end{lemma}

\begin{proof}
$q_{\varepsilon}(x)-\lambda \geq \Lambda-\lambda > 0$ for all
$x \in \Omega_{\Lambda,R}$. Therefore
\begin{align*}
  \beta(\lambda,R) 
    \geq \int_{S_{\Lambda,R}}(q_{\varepsilon}-\lambda)^2 \,\dHm_{n-1}
  \geq (\Lambda-\lambda)^2\int_{S_{\Lambda,R}} 1 \,\dHm_{n-1}.
\end{align*}
\end{proof}

\begin{corollary}
For $0 \leq \lambda < \Lambda$ and $0 \leq \rho < R \leq R_0$
\begin{equation*}
  \beta(\Lambda,\rho) \leq C \, \left(
    \frac{1}{(\Lambda-\lambda)^{\frac{4}{n-1}}(R-\rho)^2}
    + \frac{1}{(\Lambda-\lambda)^{2+\frac{4}{n-1}}} \right)
    \, \beta(\lambda,R)^{1+\frac{2}{n-1}}.
\end{equation*}
\end{corollary}

\begin{proof}
This follows directly by combining \eqref{32} and \eqref{33}.
\end{proof}

\begin{corollary}
For $\lambda_0 \geq 0$, $k=0,1,2,\ldots$ we define 
\begin{align*}
  \rho_k    &:= \frac{R_0}{2}+\frac{R_0}{2^{k}}, &
  \lambda_k &:= 2\lambda_0-\frac{\lambda_0}{2^k}, &
  J_k &:= \beta(\lambda_k,\rho_k).
\end{align*}
Then
\begin{equation}
\label{35}
  J_{k+1} \leq C(\lambda_0)  \cdot
    \left(2^{2+\frac{4}{n-1}} \right)^{k} 
    \cdot J_k^{1+\frac{2}{n-1}},
\end{equation}
where
\begin{equation*}
  C(\lambda_0):= C \cdot 
        \left(\frac{1}
               {\lambda_0^{\frac{4}{n-1}}\cdot R_0^2}
      + \frac{1}{\lambda_0^{2+\frac{4}{n-1}}}
    \right) \cdot 2^{2+\frac{4}{n-1}}.
\end{equation*}
Here $C$ is the constant from the previous corollary.
\end{corollary}

\begin{proof}
This follows immediately from the previous corollary. Just choose
$\lambda := \lambda_k$, $\Lambda:= \lambda_{k+1}$,
$\rho:= \rho_{k+1}$, $R:=\rho_k$.
Then we have $\Lambda-\lambda=\frac{\lambda_0}{2^{k+1}}$ and
$R-\rho=\frac{R_0}{2^{k+1}}$.
\end{proof}

\begin{remark} 
Reviewing the involved inequalites, we see that the constant $C(\lambda_0)$ 
in the previous corollary depends on $R_0$, the functions 
$k$, $h \in C^{\infty}(\overline{\Omega})$, as well as the smooth 
vector field $W$.
What is more important is the fact that $C(\lambda_0)$ is independent
of $\varepsilon$.
\end{remark}

\begin{proposition}
Consider $U \subset\subset \Omega$. 
Then there exists a constant $C_U$ depending on $U$ such that 
for all $0 < \varepsilon \leq 1$ the unique minimizer $u_{\varepsilon}$  
of $\AF_{\varepsilon}$ satisfies
\begin{equation*}
  \sup_U \norm{\grad u_{\varepsilon}} \leq C_U < \infty.
\end{equation*}
\end{proposition}

\begin{proof}
Without loss of generality we can restrict our attention to subsets
$U$ of the form $U:=B_{\frac{R_0}{2}}(x_0)$ with $R_0$ such that 
$B_{2R_0}(x_0) \subset \Omega$. Choosing now $\lambda_0$ large, the 
constant $C(\lambda_0)>0$ in the previous corollary can be made arbitrarily
small. Therefore, taking $\lambda_0$ large enough, we can achieve
\begin{align*}
  J_1 \leq a < 1 \quad \text{and} \quad 
  \left(2^{2+\frac{4}{n-1}} \cdot a^\frac{2}{n-1} \right)^k 
  \leq \frac{a}{C(\lambda_0)}
  \quad\text{for all } k = 0,1,2,\dotsc
\end{align*}
for an appropriate constant $0<a<1$.
Hence
\begin{equation*}
C(\lambda_0)\left(2^{2+\frac{4}{n-1}}\right)^k
  \left(a^k\right)^{1+\frac{2}{n-1}}
  \leq a^{k+1}. 
\end{equation*}
Applying \eqref{35} and proceeding by induction we obtain $J_k \leq a^k$ for
the quantity $J_k$ of the previous corollary. Now $\lim_{k \to \infty} J_k =
0$ because $a<1$ and consequently
$\beta\left(2\lambda_0, \tfrac{R_0}{2}\right) = 0$. 
\end{proof}

Summarizing our work of Sections \ref{40}, \ref{60},
and \ref{37}, Proposition \ref{39} implies

\begin{theorem}
\label{58}
The minimizing problem $\AF(u) \to \min$, 
$u \in W^{1,1}(\Omega)$, $\int_{\Omega} u \,\dvol_g = 0$ has a unique solution
\begin{equation*}
 u_0 \in W^{1,1}(\Omega) \cap C^\infty(\Omega).
\end{equation*}
\end{theorem}

In other words: Our symmetrization procedure has the desired properties. That
is, given a transvection $\tau_{\gamma}$ and a bounded subset $\hOmega \subset
\hM^n$  as described in Section \ref{4}, the symmetrization procedure
yields a unique 
set $S_{\gamma}(\hOmega)$ which has the same volume as $\hOmega$ but minimal
surface area among all sets obtained as variations of $\hOmega$ along 
the orbits of $\tau_{\gamma}$.
$S_{\gamma}(\hOmega)$ is given by
\begin{equation*}
  S_{\gamma}(\hOmega)
  = \{\tau_{\gamma}(t,\sigma(x))\;|\;
      u_0(x)-h(x) \leq t \leq u_0(x)+h(x),\; x \in \Omega\}
  \subset \hM^n. 
\end{equation*}
Using this description, we see that the boundary of $S_{\gamma}(\hOmega)$ is
smooth in points corresponding to the interior of $\Omega$.
Boundary regularity of $u_0$ on $\partial \Omega$ has not been investigated
here.

\begin{remark}
Our symmetrization construction can be carried out in any simply connected
symmetric space of nonpositive curvature. All its properties discussed above
remain valid without modifications of the proofs. In other words: An
euclidean factor does not disturb our
symmetrization constructions. Nevertheless, the case of symmetric spaces of
noncompact type is the most interesting for investigating the isoperimetric
problem. 
\end{remark}


\section{Application to the isoperimetric problem}
\label{44}



\subsection{Convexity of isoperimetric solutions}
\label{52}


Our symmetrization argument based on transvections coincides with Steiner
symmetrization in the case $\hM^n = \bbR^n$. In this case the symmetrization
procedure shows that isoperimetric solutions are convex: Suppose a
geodesic $\gamma: \bbR \to \bbR^n$ intersects a (smooth) domain $\hOmega
\subset \bbR^n$ at least twice, then the symmetrized set $S_{\gamma}(\hOmega)$
will intersect $\gamma(\bbR)$ only once. 
This is just a consequence of the fact $w \equiv 0$ if $\hM^n=\bbR^n$.
Here $w$ is the $1$--form introduced in Section \ref{40} and incorporated in
the area functional
\begin{equation*}
  \AF(u)=\int_{\Omega} \sqrt{1+k^2\norm{w+dh+du}^2}
                      +\sqrt{1+k^2\norm{w-dh+du}^2}\,\dvol_g.
\end{equation*}
In a general symmetric space $\hM^n$ of noncompact type we
would immediately get convexity of isoperimetric solutions if 
(for every geodesic $\gamma:\bbR \to \hM^n$) we had 
$w=\grad v$ for an appropriate function $v$ on the orbit space 
$M^{n-1}=\hM^n/\tau_{\gamma}$. In this case we could just set $u=-v$ to
achieve convexity.
\smallskip 

The first de Rham cohomology group of the orbit space $M^{n-1}$ is trivial,
because $M^{n-1}$ is diffeomorphic to $\bbR^{n-1}$. Therefore
$w=\grad v$ for an appropriate function $v$ if and only if $dw=0$.
Investigating for which class of symmetric spaces of noncompact type we have
$dw=0$ (for all directions of symmetrization given by geodesics $\gamma$),
it turns out that this only holds for spaces with constant sectional
curvature.
\smallskip

Suppose we are given a geodesic $\gamma:\bbR \to \hM^n$ with 
$dw \not\equiv 0$ on the corresponding orbit space
$M^{n-1}=\hM^n/\tau_{\gamma}$. Then there exists a set
$\hOmega \subset \hM^n$ which is invariant under symmetrization
with respect to $\tau_{\gamma}$ but not convex. Such an example can be
constructed as follows:

As $dw \not\equiv 0$, there exists a $2$--dimensional submanifold 
$B^2 \subset M^{n-1}$ homeomorphic to a disc such that
$\int_{B^2} dw \neq 0$. Using Stoke's theorem
\begin{equation*}
  \int_{\partial B^2} w = \int_{B^2} dw \neq 0
\end{equation*}
for the closed loop $\partial B^2$, parametrized by $s:(0,1)\to\partial B^2$.
Now we choose a function $u:M^{n-1} \to \bbR$ such that $u$ is smooth almost
everywhere, $u \circ s$ is strictly monotone increasing and ${\grad u}_{|s}$ is
tangential to $\partial B^2$ with $\norm{{\grad u}_{|s}} \equiv 1$.
Then obviously 
$\int_{\partial B^2} w(\grad u) = \int_{\partial B^2} \scp{w}{du} \neq 0$.
Choosing now a neighborhood $\Omega$ of $\partial B^2$ and a constant 
$t \in \bbR$ in an appropriate way, we can therefore achieve
\begin{align*}
  \int_{\Omega} &\sqrt{1+k^2\norm{w+d(t\cdot u)}^2}\,\dvol_g\\
   &\begin{aligned}
    &= \int_{\Omega} \sqrt{1+k^2(\norm{w}^2+2t\scp{w}{du} + t^2 \norm{du}^2)}
      \,\dvol_g\\
    &< \int_{\Omega} \sqrt{1+k^2\norm{w}^2}\,\dvol_g.
  \end{aligned}
\end{align*}
Using this it is now clear how to construct a set $\hOmega \subset \hM^n$
which is not convex but invariant under symmetrization with respect to the
transvection $\tau$. Just think about $\hOmega$ as a neighborhood of the
``lifted loop'' $s$, that is of the curve
\begin{equation*}
  \tau(t\cdot u\circ s, \sigma \circ s):(0,1) \to \hM^n.
\end{equation*}
In other words, our standard counterexample against an immediate convexity
proof by symmetrization looks like a helix winding up and
overlapping only over a very small part of the projection
$\pi(\hOmega) \subset M^{n-1}$ in the orbit space.
Of course, we expect that such a helix will not survive as a candidate for an
isoperimetric solution if we consider symmetrization with respect to another
direction, but unfortunately this is hard to control.


\subsection{Complex hyperbolic space}
\label{51}


It is well known that the boundaries of metric balls in complex hyperbolic 
spaces provide surfaces of constant mean curvature. That is, they are
critical points of the area functional with respect to volume preserving
deformations. However, for large volumes it is not known that they are
isoperimetric solutions.
Symmetrization with respect to transvections corresponds to a
special class of volume preserving deformations. Consequently metric balls
in complex hyperbolic spaces remain invariant under our symmetrization
procedure.
The fact that the area functional is convex with respect to our restricted
class of deformations provides some evidence that isoperimetric solutions
in complex hyperbolic spaces are balls and hence unique.
\smallskip

For studying our symmetrization procedure,
up to now we have only used very basic properties of the function $k$ and the
$1$--form $w$ involved in the area functional $\AF$. In fact, we did
not need much more than smoothness and $k \geq 1$.
But since we are in a symmetric space of noncompact type, these quantities
should have a lot of nice properties, remember for example Lemma
\ref{46} and \ref{47}.
This provides an interesting starting point for future research.
For intuition, we will finish this paper by explicitly computing the
$1$--form $w$ for the case of the complex hyperbolic space.

The $1$--form $w$ on the orbit space $M^{n-1}=\hM^n/\tau_{\gamma}$ is defined
by 
\begin{equation*}
  d\sigma(X) = \Hor(X)_{|\sigma} + w(X) \cdot K_{|\sigma},
\end{equation*}
where we choose $\sigma:M^{n-1} \to \hM^n$ to be the section defined in 
\ref{36} such that 
$\sigma(M^{n-1})
 =\left( \tfrac{1}{2}(\beta_{\gamma}^+ -\beta_{\gamma}^-)\right)^{-1}(\{0\})$.
Taking the scalar product with the unit normal field $\nu$ on 
$\sigma(M^{n-1})$ (where $\nu_{|\gamma(0)} = K_{|\gamma(0)}$) we obtain
\begin{align*}
  \scp{\Hor(W)}{\Hor(X)} &= \scp{W}{X} = w(X)
  = \frac{-\scp{\nu}{\Hor(X)}}{\scp{K}{\nu}}\\
  &= \frac{-1}{\scp{K}{\nu}}\left\langle 
        \nu - 
        \left\langle \nu,\tfrac{K}{\norm{K}} 
          \right\rangle \tfrac{K}{\norm{K}}
        \, ,\,
     \Hor(X) \right\rangle.
\end{align*}
Therefore
\begin{equation*}
  \Hor(W)_{|\sigma} = \frac{K_{|\sigma}}{k^2}
                     -\frac{\nu_{|\sigma}}{\scp{K}{\nu}_{|\sigma}}.
\end{equation*}
As an easy application we can compute
\begin{align*}
  \norm{W}^2 = -\frac{1}{k^2} + \frac{1}{\scp{K}{\nu}^2}
  \quad \text{and} \quad
  \sqrt{1+k^2\norm{W}^2} = \frac{k}{\abs{\scp{K}{\nu}}}.
\end{align*}

For the rest of this section we specialize to the case of 
complex hyperbolic space
$\hM^{2n}=\bbC H^n$. 
For every $X \perp K_{|\gamma(0)}$ we then
have an isometry $\varphi$ with $d\varphi(X)=X$
and $d\varphi(K_{|\gamma(0)})=-K_{|\gamma(0)}$.
These isometries can be easily obtained by direct construction.
In particular the existence of such isometries implies  
$\sigma(M^{2n-1}) = \exp(\dot{\gamma}(0)^{\perp})$.

Now we will compute the vector field $W$ explicitly. 
For this we only have to determine $\nu$ and $K$ along the geodesics 
$c:\bbR \to \hM^{2n}$ with $\dot{c}(0) \perp K_{|\gamma(0)}$.
As $\hM^{2n}$ is the complex hyperbolic space, the operator
$R(\cdot, \dot{c}(t))\dot{c}(t)$ has eigenvalues $0$, $-1$ and $-4$. More
precisely, along $c$ we may choose $2n$ orthonormal parallel vector fields 
$X_1, JX_1=\dot{c}, \dotsc, X_n, JX_n$, where $J$ denotes the almost complex
structure, such that  
\begin{align*}
  R(X_1(t),\dot{c}(t))\dot{c}(t) &= - 4 X_1(t), \\  
  R(JX_1(t),\dot{c}(t))\dot{c}(t) &= 0, \\
  R(X_i(t),\dot{c}(t))\dot{c}(t) &= - X_i(t), && \text{for } i\geq 2\\
  R(JX_i(t),\dot{c}(t))\dot{c}(t) &= - JX_i(t), && \text{for } i\geq 2.
\end{align*}
Furthermore, we can assume 
$\dot{\gamma}(0) = \cos(\vartheta) \cdot X_1(0)+\sin(\vartheta)\cdot X_2(0)$
for an appropriate $\vartheta \in [0,2\pi)$ without restriction.
Since $K$ is a Killing field corresponding to a transvection, 
$K(t):=K_{|c(t)}$ is a Jacobi field with initial data
$K(0)=\dot{\gamma}(0)$ and $K'(0)=0$. This implies
\begin{equation*}
  K(t) = \cos(\vartheta) \cosh(2t) X_1(t) 
          + \sin(\vartheta) \cosh(t) X_2(t).
\end{equation*}
For calculating the normal field $\nu$ along $c$, we consider 
curves of the form 
$Y: (-\varepsilon,\varepsilon) \to T_{\gamma(0)}\sigma(M^{2n-1})$
with $\norm{Y(s)} \equiv 1$. Applying the usual Jacobi field techniques
to the ``radial'' geodesic variations
$V(t,s):=\exp_{\gamma(0)}(t\cdot Y(s))$, it turns out that 
$T_{c(t)}\sigma(M^{2n-1})$ is spanned by the Jacobi fields
\begin{align*}
  & JX_1(t) =  \dot{c}(t), \\
  & \sinh(t) \cdot X_i(t) \phantom{J}\qquad \text{ for } i=3,\dotsc,n,\\
  & \sinh(t)  \cdot JX_i(t) \qquad \text{ for } i=2,\dotsc,n, \\
  & \sin(\vartheta) \cdot \tfrac{1}{2}\sinh(2t) \cdot X_1(t)
    - \cos(\vartheta) \sinh(t) \cdot X_2(t).
\end{align*}
As $\nu(t):=\nu_{|c(t)}$ has to be perpendicular to these vector
fields
\begin{equation*}
  \nu(t) = \frac{\cos(\vartheta)\cdot X_1(t)
                 +\sin(\vartheta) \cosh(t) \cdot X_2(t)}
              {\sqrt{\sin^2(\vartheta) \cdot \cosh^2(t) +\cos^2(\vartheta)}}.
\end{equation*}
Combining these results we obtain
\begin{equation}
\label{45}
\begin{split}
  \Hor(W)_{|c(t)} =
  \begin{aligned}[t]
    & \frac{\cos(\vartheta)\cosh(2t)\cdot X_1(t) 
            + \sin(\vartheta)\cosh(t) \cdot X_2(t)}
           {\cos^2(\vartheta)\cosh^2(2t)+\sin^2(\vartheta)\cosh^2(t)}\\[1ex]
    & - \frac{\cos(\vartheta) \cdot X_1(t)
              +\sin(\vartheta) \cosh(t) \cdot X_2(t)}
	   {\cos^2(\vartheta)\cosh(2t)+\sin^2(\vartheta)\cosh^2(t)}
  \end{aligned}
\end{split}
\end{equation}
Observe that 
$\cos(\vartheta) = \scp{\dot{\gamma}(0)}{X_1(0)}
 = \scp{J\dot{\gamma}(0)}{JX_1(0)}
 = \scp{J\dot{\gamma}(0)}{\dot{c}(0)}$. In other words,
$\vartheta$ is the angle between $J\dot{\gamma}(0)$ and $\dot{c}(0)$.
Using \eqref{45} we can immediately derive the following qualitative
properties of $\Hor(W)_{|c(t)}$:
\begin{itemize}
\item $\Hor(W)_{|c(t)} \perp \dot{c}(t)$.
\item $\Hor(W)_{|c(t)} \to 0$ for $t \to \infty$.
\item $\Hor(W)_{|c(t)} = 0$ for $\vartheta = 0$ and $\vartheta =
  \frac{\pi}{2}$. 
\end{itemize}

\begin{remark}
The $1$--form $w$ that appears in our symmetrization construction clearly has
some interesting special properties with respect to the geometry of the
symmetric space.
This information might provide further insight into the shape of isoperimetric
domains.
There also might be interesting connections to stability of isoperimetric
domains $\hOmega$. After all, Killing fields on $\hM^n$ induce elements in the
kernel of the Jacobi operator of $\partial \hOmega$.
\end{remark}


\end{document}